%
%

\magnification=1200

\font\titfont=cmr10 at 12 pt

\font\headfont=cmr10 at 12 pt



\overfullrule=0in

\def\boxit#1{\hbox{\vrule
 \vtop{%
  \vbox{\hrule\kern 2pt %
     \hbox{\kern 2pt #1\kern 2pt}}%
   \kern 2pt \hrule }%
  \vrule}}

  \def\harr#1#2{\ \smash{\mathop{\hbox to .3in{\rightarrowfill}}\limits^{\scriptstyle#1}_{\scriptstyle#2}}\ }

\def\ss{\subset}
\def\half{\hbox{${1\over 2}$}}
\def\smfrac#1#2{\hbox{${#1\over #2}$}}
\def\oa#1{\overrightarrow #1}

\def\dist{{\rm dist}}

\def\log{{\rm log}}
\def\Hess{{\rm Hess}}

\def\tr{{\rm tr}}
\def\max{{\rm max}}

\def\span{{\rm span\,}}

\def\det{{\rm det}}

\def\Sym{{\rm Sym}^2}

\def\rn{\bbr^n}

\def\Int{{\rm Int}}

\def\Symn{{\Sym(\rn)}}
\def\SymN{{\Sym(\bbr^N)}}

\def\SA{{\rm SA}}
\def\ct{{$C^2$}}

\def\Theorem#1{\medskip\noindent {\bf THEOREM \bf #1.}}
\def\Prop#1{\medskip\noindent {\bf PROPOSITION #1.}}
\def\Cor#1{\medskip\noindent {\bf COROLLARY #1.}}
\def\Lemma#1{\medskip\noindent {\bf Lemma #1.}}
\def\Remark#1{\medskip\noindent {\bf Remark #1.}}
\def\Note#1{\medskip\noindent {\bf Note #1.}}
\def\Def#1{\medskip\noindent {\bf Definition #1.}}

\def\Ex#1{\medskip\noindent {\bf Example \bf    #1.}}

\def\pf{\medskip\noindent {\bf Proof.}\ }
\def\qed{\hfill  $\vrule width5pt height5pt depth0pt$}
\def\n{\nabla}

   \def\cp{{\cal P}}   
\def\cw{{\cal W}}

\def\cp{{\cal P}}
\def\cf{{\cal F}}

\def\vf{\varphi}

\def\wt{\widetilde}
\def\wh{\widehat}

\def\and{\qquad {\rm and} \qquad}

\def\bbr{{\bf R}}\def\bbh{{\bf H}}
\def\bbc{{\bf C}}

\def\bbz{{\bf Z}}

\def\a{\alpha}
\def\b{\beta}
\def\d{\delta}
\def\e{\epsilon}
\def\f{\phi}

\def\l{\lambda}
\def\o{\omega}

\def\s{\sigma}
\def\x{\xi}

\def\D{\Delta}
\def\L{\Lambda}
\def\G{\Gamma}
\def\O{\Omega}

\def\psh{plurisubharmonic }

\def\lloc{L^1_{\rm loc}}

\def\bo{\partial \Omega}
\def\Ob{\overline{\O}}

\def\MA{MA}

\def\Symn{\Sym(\rn)}
 
\def\USC{{\rm USC}}
\def\fa{{\rm\ \  for\ all\ }}

\def\cpt{\wt{\cp}}
\def\ft{\wt F}
\def\ob{\overline{\O}}

\def\K{{\rm K}}

\def\M{{\bf M}}
\def\N#1{C_{#1}}
\def\ds{Dirichlet set }
\def\dir{Dirichlet }
\def\Fa{{\oa F}}

 \def\LAG{{\rm LAG}}
 \def\ISO{{\rm ISO_p}}

\def\AA{1}
\def\BB{2}
\def\CC{3}
\def\DD{4}
\def\EE{5}
\def\FF{6}
\def\GG{7}
\def\GGG{8}
\def\HH{9}
\def\II{10}

\def\CD{C}


\ 
\vskip .3in

\centerline{\titfont DIRICHLET DUALITY  }
\medskip

\centerline{\titfont  and the  }\medskip

\centerline{\titfont    NONLINEAR DIRICHLET PROBLEM}
\vskip .3in

\centerline{\titfont F. Reese Harvey and H. Blaine Lawson, Jr.$^*$}
\vglue .9cm
\smallbreak\footnote{}{ $ {} \sp{ *}{\rm Partially}$  supported by
the N.S.F. }

\vskip .5in
\centerline{\bf ABSTRACT} \medskip
  \font\abstractfont=cmr10 at 10 pt
We study the \dir problem for fully nonlinear, degenerate elliptic 
equations of the form ${\bf F}(\Hess \, u)=0$ on a smoothly bounded domain
$\O\ss\ss\rn$. In our approach the equation is replaced by a subset $F\ss\Symn$ of the
symmetric $n\times n$ matrices with $\partial F \subseteq \{  {\bf F}=0\}$.
We establish  the existence and uniqueness of continuous solutions
under  an explicit  geometric ``$F$-convexity'' assumption  on the
boundary  $\bo$. 
The topological structure of $F$-convex domains is also studied and 
a theorem of Andreotti-Frankel type is proved for them.
Two key ingredients in the analysis are the use of 
subaffine functions and \dir duality, both introduced here.
Associated to $F$  is a \dir dual set $\ft$  which gives a dual \dir
problem.  This pairing is a true duality in that
the dual of $\ft$ is $F$ and in the analysis the roles of $F$ and $\ft$ are interchangeable. 
The duality also clarifies many features of the problem including the appropriate conditions
on the boundary.
Many interesting
examples are covered by these results including:  All branches of the homogeneous 
Monge-Amp\`ere equation over $\bbr$, $\bbc$ and $\bbh$; equations appearing naturally in
calibrated geometry, Lagrangian geometry and $p$-convex riemannian geometry, and
all branches of the Special Lagrangian potential equation.
{{\parindent= .43in
\narrower\abstractfont \noindent

}}

\vfill\eject

\centerline{\bf TABLE OF CONTENTS} \bigskip

{{\parindent= .1in\narrower\abstractfont \noindent

\qquad \AA. Introduction.\smallskip

\qquad \BB.     The Maximum Principle and Subaffine Functions.    \smallskip

\qquad \CC.     Dirichlet Sets -- The Maximum of Two Functions.  \smallskip

\qquad \DD. Dirichlet Duality.\smallskip

\qquad \EE.     Boundary Convexity.    \smallskip

\qquad \FF.     The Dirichlet Problem.  \smallskip

\qquad \GG.     Quasiconvex Functions.  \smallskip

\qquad \GGG.     Sup-Convolution  Approximation.  \smallskip

\qquad \HH.     Topological Restrictions on Domains with Strictly $\Fa$-Convex Boundaries.  \smallskip

\qquad \II.     Examples of Dirichlet Sets.  \smallskip

}}

\vskip .3in

{{\parindent= .3in\narrower

\noindent
{\bf Appendices: }\medskip

A. Dirichlet Sets which can be Defined using Fewer of the Variables.
\smallskip

B. A Distributional Definition of Type $F$ for Convex Dirichlet Sets in $\rn$.
\smallskip

}}

\vfill\eject

\noindent{\headfont \AA. Introduction.} \smallskip

The point of this paper is to study the \dir problem for certain fully nonlinear, degenerate
elliptic, second order differential equations which appear naturally in geometry.  The class of problems 
we consider has a rich structure and covers a wide variety of interesting cases.
To be more specific,  we suppose that $\O$ is a  bounded domain  in $\rn$ with smooth boundary
and that the  nonlinear operator ${\bf F}$  depends only on the second derivatives of the
unknown function.  We then consider the homogeneous \dir problem: to show
$$
{\rm Given\ }  \vf \in C(\bo), \ \ \exists \, ! \ u\in C(\overline \O)\ \ {\rm with\ \ }{\bf F}(\Hess \, u)=0\  {\rm on\ } \O \ {\rm \  and\ }\ u\bigr|_{\bo} = \vf.
\leqno{(DP)}
$$
To our surprise, uniqueness does not seem to be included in the celebrated theory of 
viscosity solutions unless ${\bf F}$ is either uniformly elliptic or proper with respect
to the variable $u$.  Moreover, a local geometric condition
on $\bo$ needed for existence only seems to be available in certain cases (cf. the inspiring
paper [CNS]). We shall give answers to these two questions.  

We  take a  geometric approach to the equation (in the spirit of Krylov [Kr]) 
which eliminates the operator ${\bf F}$ and replaces it with a closed subset $F$ of 
the space $\Symn$ of real symmetric $n\times n$ matrices, with the property  that $\partial F$ is 
contained in $\{ {\bf F}=0\}$. 
In this approach we formulate the notion of solution as  a  {\sl dual notion}. 
Although the fact is not needed in this paper,   we show at the end of Section 4 that our solutions are  the standard viscosity solutions.  We feel our duality makes all the basic properties and the comparison theorems more transparent.  Furthermore,  this duality is a true duality in that every equation has a well defined dual equation, and their roles are interchangeable in the theory.

The geometric approach to the problem  also leads naturally to a pointwise convexity condition
on the boundary $\bo$ needed for the existence question.  This condition generalizes 
the usual convexity and pseudoconvexity for the classical Monge-Ampere equation in the
real and complex case, as well as the $\phi$-convexity introduced for domains in a calibrated
manifold $(M,\phi)$ in [HL$_2$].

Interestingly, this   convexity condition for $\bo$  gives explicit restrictions on the topology of the domain $\O$.  In particular, there is an integer $D$, depending only on the  subset $F$,  such  that if $\bo$ satisfies the  convexity condition at each point, then $\O$ has 
the homotopy type of a CW-complex of dimension $\leq D$.

An important aspect of this theory is that it applies to a wide spectrum of interesting cases.
For example, suppose $K = \bbr, \bbc$ or $\bbh$ (the real, complex and quaternionic number fields
respectively), and consider $K^n=\bbr^N$ where $N=n, 2n$ or $4n$.  Every real symmetric 
$N\times N$-matrix $A$ has a $K$-hermitian symmetric part $A_K$ with eigenvalues
$\l_1\leq \cdots \leq \l_n$.  
%
%
The associated determinant $\det_K A_K = \l_1\cdots\l_n$ is a polynomial of degree $n$  in the entries of $A$ and there is an associated Monge-Amp\`ere equation
$$
\det_K  \{\Hess\, u\}_K\ =\ 0.
$$
Solutions to the \dir problem for this equation are understood in the case where $\{\Hess\,u\}_K\geq0$, i.e., $\l_1=0$
  (see, for example, [Alex$_2$],  [Br],  [BT],  [Al]).  However, our theory
gives unique solutions of (DP) for the other branches of the equation, namely
$$
\l_q\ =\ 0
$$
for any fixed $q$.  This important result is due to  Hunt and Murray [HM] and Slodkowski [S]
in the complex case.  The  work of Slodkowski was an inspiration for this paper.   
His result on the largest eigenvalue of a convex function is the 
deepest ingredient in our uniqueness proof.  

We note incidentally that the problem dual to $\l_q=0$ in our sense is $\l_{n-q+1}=0$.

One can also treat the equation
$$
\l_p+\l_{p+1}+\cdots+\l_{p+q}\ =\ 0
$$
for fixed $p$ and $q$ by these methods.
 
 A large and  important class of examples are those which are {\sl geometrically defined}.
In particular every calibration on $\rn$ gives rise to an equation of our type.
More details are given just below.

Yet another interesting case is the equation
$$
{\rm Im} \left\{ e^{i\theta} \det (I+i\Hess\, u)\right\} \ =\ 0,
$$
(for fixed $\theta$) which governs the potential functions in the theory of Special Lagrangian
submanifolds. The locus of this equation, considered as a subset of $\Symn$, has 
$n$ distinct connected components or {\sl branches}, unless $\theta =\pi/2$ when $n$ is odd  or 
$\theta=0$ when $n$  is even. In these exceptional cases there are $n-1$ branches.
For the two outermost branches and with $\theta=0$, the \dir problem was treated in depth in
[CNS].  Furthermore, they conjectured that there exist the same number of solutions as there are branches.  Our results  show that indeed the \dir problem is uniquely solvable in 
continuous functions for every branch and for every   $\theta$   in each dimension.  
In particular the $n$ (or $n-1$) distinct solutions for a given boundary function exist and are uniquely determined by the   distinct branches.
They are also nested, i.e., $u_1\leq u_2\leq \cdots$.

Our general set-up here is the following.  
We start with a given closed subset $F$ of the space real symmetric matrices $\Symn$.
We are interested in formulating and  solving the \dir problem for the equation
$$
\Hess_x  u \in \partial F \fa x\in \O.
\eqno{(\AA.1)}
$$
using the functions of ``type $F$'', i.e., which satisfy 
$$
\Hess_x u \in   F \fa x.
\eqno{(\AA.2)}
$$
{\sl Apriori} these conditions make sense only for $C^2$ functions $u$.  
We shall extend the notion to functions which are only upper semicontinuous.

  This extension requires two ingredients. First we introduce the class of 
  {\bf subaffine functions}.  These are upper semicontinuous functions $u$ defined locally by 
  the condition:
  \medskip

{\sl For each affine function $a$, if $u\leq a$ on the boundary of a ball $B$, then $u\leq a$ on $B$.}

\medskip

These locally subaffine functions are globally subaffine and hence satisfy the maximum
principle on any compact set.  A $C^2$-function is subaffine if and only if $\Hess\,u$ has at
least one eigenvalue $\geq0$ at each point.

The second  step is to consider the {\bf \dir dual set}
$$
\ft \ \equiv\ -(\sim \Int F).
\eqno{(\AA.3)}
$$
and  define an upper semicontinuous function $u$ to be {\bf of type} $F$ if 
\bigskip
\centerline{
$u+v$  is subaffine for all
$C^2$-functions $v$ of type $\ft$.
}
\bigskip
\noindent
 In other words,  $u\in \USC$ is type $F$ if   for any  ``test function''   $v\in C^2$ of dual type $\ft$,
the sum $u+v$ satisfies the maximum principle. 
 
 Note that   $\wt{\ft}=F$, and so our condition above has an inherent symmetry.

The key requirement on $F$ for solving the \dir problem
 is that the maximum of two functions of type $F$ be again of type $F$.  
 This is, in effect, equivalent to the following positivity condition on our set.  We say that  $F$  is  
 a {\bf \dir set} if it satisfies the condition
$$
F+\cp \ss F.
\eqno{(\AA.4)}
$$
where 
$$
\cp\  =\ \{A\in\Symn : A\geq0\}
$$
 is the subset of non-negative matrices.  This condition corresponds to degenerate ellipticity in modern fully nonlinear theory.

The simplest case, where $F=\cp$, is classical.  Here the functions of type $\cp$ are the convex functions, and strict $\cp$-convexity of the boundary is the conventional notion.

In the dual case  where  $F =\cpt =  \{A\in\Symn : A\not < 0\}$ we shall prove that  an upper semicontinuous function $u$ is of type $\cpt$ if and only if it is subaffine.

It is easy to see that $F$ is a \dir set if and only if $\ft$ is a \dir set. 

\dir sets can be quite general in structure. Translates, unions (when closed)
 and intersections  of \dir sets are \dir sets. However there are quite interesting ones coming from geometry as follows.  Let $G(p,\rn)$ denote the 
grassmannian of $p$-planes in $\rn$ and fix any compact subset $G\ss G(p,\rn)$. Let
$$
\cp(G)\ =\ \{A\in\Symn : {\rm trace}(A\bigr|_\x)\geq 0 \fa \x\in G\}
$$
Then $\cp(G)$ is a \dir set with \dir dual
$$
\cpt(G)\ =\ \{A\in\Symn : {\rm trace}(A\bigr|_\x)\geq 0 \ \ {\rm for\ some\ } \x\in G\}
$$
The $C^2$-functions of type $\cp(G)$ are characterized by being subharmonic on all $G$-planes.
In fact they are subharmonic on all minimal $G$-submanifolds (those whose tangent planes are
all $G$-planes). Every calibration $\phi$ gives a set $G=G(\phi)$ of this type where $G$-submanifolds are automatically minimal. As special cases one considers
the complex and quaternionic grassmannians.  Another interesting case, not coming from a calibration,  is given by the set 
$G=LAG$  of all Lagrangian $n$-planes in $\bbc^n$.

Further interesting examples arise from restriction.  If $F_W\ss\Sym(W)$ is a \dir set, where 
$W\ss\rn$ is a linear subspace, then $F=\{A\in\Symn : A\bigr|_W\in F_W\}$ is also a \dir set. 
Since arbitrary intersections of \dir sets are \dir sets, this yields a new  \dir set for each family 
of \dir sets on subspaces of $\rn$

In addition, many of the  interesting examples can be introduced in terms of
 G\aa rding polynomials on $\Symn$ with the identity $I$ a hyperbolic direction.
 These in turn can generate more examples by taking directional derivatives in the direction $I$.

The very general nature of \dir sets complicates the question: ``What geometric conditions on $\bo$ are necessary to solve the \dir problem for a given $F$?''  \ Associated to each $F$ is an asymptotic 
cone or ``ray set'' $\Fa$. This is a closed cone with vertex at the origin  and
consisting essentially of the rays which   lie inside $F$ after some point. 

Suppose now that  $\O\ss\ss\rn$ is a domain with smooth boundary.  Denote by $II$ the second
second fundamental form of the boundary with respect to the inward-pointing unit normal. 
Then $\bo$ is said to be  {\bf strictly $\Fa$-convex} at $x\in \bo$ if
$$
II_x  \ =\ B\bigr|_T  \ \ {\rm for\ some\ } B\in \Int\Fa.
$$
where $T=T_x(\bo)$.  This is equivalent to the condition that 
$II_x +tP_n\in \Int\Fa$ for all $t\geq $ some $t_0$ where $P_n$ is projection onto the 
line normal to $\bo$ at $x$.

By a {\sl global defining function} for $\bo$ we mean   a function $\rho\in C^\infty(\ob)$ with $\rho<0$ on $\O$ and with  $\rho=0$ and $\n \rho\neq 0$ on $\bo$. 
We prove the following result.

\Theorem{\EE.12}  {\sl Suppose $F$ is a Dirichlet set.
 If the boundary $\bo$ is strictly $\Fa$-convex at each point,
then there exists  a global defining function $\rho \in C^\infty(\Ob)$
 for $\bo$  which is strictly of type $\Fa$
on $\Ob$.  Moreover, }
$$
 \exists\ \e>0 \  {\rm and }  \ R>0 \ {\rm such\  that\ }\  C\left(  \rho - \e \half |x|^2\right)\ \in \   F(\Ob) \fa C\geq R
 $$

We are now in a position to discuss the main theorem.  A function $u$ on a domain $\O$ 
is said to be {\sl $F$-\dir } if  $u$ is of type $F$ and $-u$ is of type $\ft$. Such a function $u$ is automatically continuous, and at any point $x$ where $u$ is $C^2$, it satisfies the condition (\AA.2) above.

\Theorem{\FF.2. (The Dirichlet Problem)}  {\sl  Let $\O$ be a bounded domain in $\rn$
with  smooth boundary $\bo$, and let $F$ be a Dirichlet set.
Suppose that $\bo$ is both $\Fa$ and $\oa \ft$ strictly convex at each point. Then  for each $\vf\in C(\bo)$, there exists a unique  $u\in  C(\overline\O)$ which is an $F$-Dirichlet function  on $\O$ and equals $\vf$ on $\bo$.
}

\medskip
\noindent
{\bf Note.}  The requirement of both $\Fa$-convexity and $\oa \ft$-convexity for $\bo$
is necessary.  In fact this explains the restriction $2q<n$ which appears in the work of
Hunt and Murray [HM].

\medskip
\noindent
{\bf Note.}  Well known uniqueness results (cf. [J], [I], [IL], [SO$_1$] for example) require either 
uniform ellipticity  or properness of the equation with respect to the variable $u$.
See [CIL], and [SO$_*$] for a fuller account and references.

\medskip

The uniqueness part of Theorem \FF.2 is immediate  from the following comparison result
and the maximum principle for subaffine functions.
For an open set $X\ss\rn$, let $F(X)$ denote  the set of (u.s.c.) functions of type $F$ on $X$
and let $\SA(X)$ denote the subaffine functions on $X$.

\Theorem{\FF.5. (The Subaffine Theorem)}  
{\sl  Assume that $F$ is a Dirichlet set.

\noindent
 If $u\in F(X)$ and $v\in\ft(X)$, then $u+v \in\SA(X)$.}

\medskip

The proof of this result    is given in  Sections \GG \ and \GGG.
In Section \GG \ we use a breakthrough technique of 
Slodkowski to prove the Subaffine Theorem when $u$ and $v$ are quasi-convex.
Slodkowski's  work enables one to pass from an  estimate which holds 
almost everywhere to one which holds at all points, and can therefore be used to
establish the maximum principle.
Then, in Section \GGG, \  sup-convolution is used to approximate arbitrary
$u$ and $v$, of type $F$ and $\ft$ respectively, by quasi-convex functions of the same type.

\medskip
\noindent
{\bf Remark.}   We note that if $F_1\ss F_2$ are \dir sets, and if $u_1, u_2$ are the corresponding
solutions to the \dir problem above (for the same  boundary function $\vf$), then
$$
u_1\ \leq\ u_2 \quad{\rm on}\ \ob.
\eqno{(\AA.5)}
$$
Thus the entire lattice of \dir sets, ordered by inclusion, maps in an order preserving way
to the set of solutions.
If, for example, one restricts to \dir sets which are cones with vertex at the origin, then
our ordered family has an initial object $\cp$ and final object $\cpt$. 
For any continuous function given on the boundary of a convex domain
we obtain a huge family of solutions all lying above the convex solution and below the 
concave one.  They serve as ``barriers'' for each other in that (\AA.5) holds whenever
$F_1\ss F_2$.  Of course, somewhere in there lies the harmonic solution corresponding to
$F_{\rm harm} = \{A: \tr A\geq 0\}$. Even in two variables it is interesting to contemplate this family.
Within it, for example, lie the \dir sets $F=\{A : a_{11}\geq 0\}$ and
$F=\{A : a_{11}\geq 0\ {\rm and\ } a_{22}\geq 0\}$
whose associated \dir functions are weak solutions of $u_{xx}=0$ and $u_{xx}u_{yy}=0$
respectively. 

\Remark{}  The case $F=\{A : a_{11}\geq 0\}$, corresponding to $u_{xx}=0$, demonstrates
 the utter lack of regularity (beyond continuity) for general solutions obtained here. 
 The $F$-convexity required for a domain $\O\ss\bbr^2$ is that it be horizontal-slice convex 
 (i.e., horizontal slices are connected) and the unique solution for a given boundary function $\vf$ is the linear interpolation
 on these slices.

\medskip
The paper is organized as follows.

In Section 2 we introduce the notion of a subaffine function.  This is a class  of functions which satisfy the maximum principle and are determined by local properties.

In Section 3 the ``positivity condition'' $F+\cp\ss F$ is discussed in some detail.  For convenience, and to avoid the overused word {\sl elliptic}, these sets are called \dir sets.  This is exactly the 
natural condition to ensure that $u,v\in F(X) \ \Rightarrow\ \max\{u,v\}\in F(X)$.

In Section 4 the dual set $\ft$  is investigated.  This duality clarifies our weak definition of
type $F$ and leads to a natural discussion of uniqueness via the Subaffine Theorem.

In Section 5 the associated ray set $\Fa$ is introduced, $\Fa$-convexity of the boundary is discussed,  and Theorem \EE.12 is proved.

In Sections 6, 7 and 8 the  \dir problem is solved.  Existence follows from the Perron method  and the classical 
``barrier'' argument, combined with a regularity argument of Walsh.  Uniqueness is reduced to the Subaffine Theorem, which is proved in Sections 7 and 8. 

In Section 9 we show that the natural domains $\O$ for which the $F$-\dir problem can be solved
are topologically restricted.  If $D$ is the ``free dimension of $F$'', then $\O$ has homotopy dimension $D$ and   $H_k(\bo;\,\bbz) \cong H_k(\O;\,\bbz) \fa k<n-D-1$.

In Section 10 we discuss numerous examples of  \dir sets, as well as general principles for constructing them.  This section shows that there are many interesting applications of the main results.

In Appendix A we show that for  \dir sets $F$ which can be defined using 
fewer of the variables in $\rn$ ( i.e., in terms of a \dir set $F_0$ associated
 to a proper subspace  $\bbr^p\ss\rn$) that  an u.s.c.  function $u$ is of type $F$ 
 if and only if the restriction of $u$  to each horizontal $\bbr^p$ is of type $F_0$.
 
 In Appendix B  a distributional definition of type $F$ is given when $F$ is convex.

\Note{}  Throughout this paper $X$ will denote an open connected subset of $\rn$.

\vfill\eject


\noindent{\headfont \BB. The Maximum Principle and Subaffine Functions.} \smallskip

 For a discussion of the maximum principle it is natural to consider the space $\USC(X)$
 of upper semi-continuous functions on $X$ with values in $[-\infty,\infty)$.
 A  function $u\in \USC(X)$  satisfies the {\sl maximum principle} if for each 
 compact subset $K\ss X$
$$
\sup_K u\ \leq \ \sup_{\partial K} u.
\eqno{(\BB.1)}
$$
 A function $u$ may locally satisfy the maximum principle
  without satisfying the maximum principle on all of $X$.
(Consider, for example, a function $u$ on $\bbr$ with compact support, $0\leq u\leq 1$,
$u\equiv 1$ on a neighborhood of the origin and otherwise monotone.)  However, this situation is easily remedied.  First note that (\BB.1) is equivalent 
to the condition that:
$$
u\ \leq \  c \ \ {\rm on\ } \partial K \qquad\Rightarrow\qquad 
u\ \leq \  c \ \ {\rm on\ } K  \qquad
 {\rm for\ all\ constants\ } c,
\eqno{(\BB.1)'}
$$
 i.e., $u$ is {\sl sub-constants}.   Replacing the constant functions by the affine functions, consider the condition:
 $$
u\ \leq \  a \ \ {\rm on\ } \partial K \qquad\Rightarrow\qquad 
u\ \leq \  a \ \ {\rm on\ } K \qquad
 {\rm for\ all\ affine\ functions\ } a
\eqno{(\BB.2)}
$$
\Def{\BB.1}   A function $u\in \USC(X)$ satisfying (\BB.2)
 for all compact subsets $K\ss X$  will be called {\sl subaffine on}  $X$. 
Let $\SA(X)$ denote the space of all $u\in \USC(X)$ that are locally subaffine on $X$,
i.e., for all $x\in X$ there exists a neighborhood $B$ of $x$ with $u\bigr|_B$ subaffine on $B$.
\medskip

Note that if $u$ is subaffine on $X$, then the restriction to any open subset is also subaffine.

\Lemma{\BB.2}  {\sl If $u$ is locally subaffine on $X$, then $u$ is subaffine on $X$.

\noindent
Moreover, $u\notin \SA(X)$   if and only if }
$$\eqalign{
{\sl there\ exist\ \ }  x_0\in X, \ \ &a\ {\sl affine, \ and\ } \e>0 \ {\sl such\ that\ } \cr
 (u-a)(x)\ \ &\leq\  -\e|x-x_0|^2\ \ {\sl near\ } x_0, \ \ {\sl and}\cr
 (u-a)(x_0)\   &=  \  0\  \cr
 }
\eqno{(\BB.3)}
$$
\pf
Subaffine implies locally subaffine, which implies (\BB.3) is impossible.  
Hence, it remains to show that if (\BB.3) is false, then $u$ is subaffine, or
equivalently,
 if $u$ is not subaffine on $X$, then (\BB.3) is true.
If $u$ is not subaffine on $X$, then for some compact set $K\ss X$ and some affine function
$b$, the difference $w=u-b$ has a  strict interior maximum point for $K$.  For $\e>0$ sufficiently small,
the same is true for $w = u+ \e |x|^2-b$.    Choose a maximum point $x_0\in\Int K$ for $w$
and let $M=w(x_0)$ denote the maximum value on $K$.    
Then $u+ \e |x|^2-b-M \leq 0$ on $K$ and equals zero at $x_0$.  Since $\e|x|^2$ and $\e|x-x_0|^2$
differ by an affine function, this proves that there is an affine function $a$ such that
$u+\e|x-x_0|^2-a\leq0$ on $K$ and is equal to zero at $x_0$, i.e., 
(\BB.3) is true. \qed

\Prop{\BB.3. (Maximum Principle)}
{\sl  Suppose $K\ss\rn$ is compact and $u\in  \USC(K)$.  If 
$u\in\SA(\Int K)$, then}
$$
\sup_K u\ \leq\ \sup_{\partial K} u.
$$

\pf  Exhaust $\Int K$ by compact sets $K_\e$. Since $u\in \SA(\Int K)$, 
$\sup_{K_\e} u\leq \sup_{\partial K_\e} u$. Since $u\in \USC(K)$, 
each $U_\d = \{x\in K: u(x) < \sup_{\partial K} u + \d\}$,
for $\d>0$, is an open neighborhood of $\partial K$ in $K$.
Therefore, there exits $\e>0$ with $\partial K_\e\ss U_\d$ which implies that 
$\sup_{\partial K_\e} u \leq \sup_{\partial K} u+\d$.  \qed

\medskip

For functions which are $C^2$ (twice continuously differentiable), the subaffine condition is a condition on the hessian of $u$ at each point.

\Prop{\BB.4}  {\sl Suppose $u\in C^2(X)$.  Then }
$$
u\in\SA(X) \ \ \ \iff \ \ \   \Hess_x u \ \ {\sl has \ at \ least \ one \ 
eigenvalue\ }\geq0 \ {\sl at\ each \ point\ }x\in X.
$$
\pf  Suppose $ \Hess_{x_0} u <0$ (negative definite) at  some point $x_0\in X$.
Then the Taylor expansion of $u$ about $x_0$ implies (\BB.3).
Therefore, since $u(x_0)=0$, $u\notin  \SA(X)$.

Conversely, if $u\notin \SA(X)$, then (\BB.3) is valid for some point $x_0\in X$ which implies
that $\Hess_{x_0} u + \e I \leq 0$.     So  $\Hess_x u <0$ is negative definite.
\qed

\medskip

\Ex{(n=1)}  Suppose $I$ is an open interval in $\bbr$.  Then
$$
u\in \SA(I) \quad\iff\quad {\rm either \ \ } u\in {\rm Convex}(I)\ {\rm or\ } u\equiv-\infty.
$$
\pf 
Suppose $u\in\SA(I)$ equals $-\infty$ at one point $\a\in I$ but $u$ is finite at another point
$\b\in I$ with $\a<\b$. Choose $a$ to be the affine function with $a(\a) = -N$ and $a(\b) = u(\b)$.
By (\BB.2),  we have $u\leq a$ on $[\a,\b]$, which implies (by letting $N\to\infty$) that
$u\equiv -\infty$ on $[\a,\b)$.  The case $\b<\a$ is identical. 
Hence $u$ is either $\equiv -\infty$ or it is finite-valued 
on all of $I$ (and therefore convex).  
The converse is immediate.\qed

\medskip

Next we give a characterization of   subaffine  functions and convex functions  which 
is the   prototype of   Dirichlet duality. 

\Prop{\BB.5}  {\sl  Suppose $v\in\USC(X)$.  Then $v\in \SA(X) \iff u+v\in \SA(X)$ for all 
$u\in {\rm Convex}(X)$.}

\pf Since $u=0$ is convex, we need only prove that if $u\in {\rm Convex}(X)$
and $v\in\SA(X)$, then $u+v\in\SA(X)$.
Equivalently,
we must show that if $v\in\SA(X)$, then:
$$
{\rm For\ all\ }  w\in {\rm Concave}(X), \ \ 
v\ \leq\  w \ \ {\rm on \ }\partial B \ \ \Rightarrow\ \  v \ \leq \ w \ \ {\rm on \ }  B
\eqno{(\BB.4)}
$$
for an arbitrary closed ball $B$ contained in $X$.  That is,
$$
v \ \ {\rm is\ subaffine\ } \qquad \Rightarrow\qquad v \ \ {\rm is\ ``subconcave''\ }
\eqno{(\BB.5)}
$$
To prove (\BB.4), choose $a$ affine with $w\leq a$ on $B$.
Then $v\leq w \leq a$ on $\partial B$ implies $v\leq a$ on $B$ since 
$v$ is subaffine. Now any concave function $w$ is the infimum over the 
family of affine functions $a$ with $w \leq a$. (To see this, apply the finite-dimensional
Hahn-Banach Theorem to the graph of $w$.)
It follows that  that $v\leq w$ on $B$.
\qed
\medskip

Let $\overline{\rm C}{\rm onvex}(X)$ denote the set of functions on $X$ that locally are 
either convex or $\equiv -\infty$.  It is easy to see by the example above, that: 
$$
u\in \overline{\rm C}{\rm onvex}(X)\quad  \iff \quad {\rm  the\ restriction\ of\ }  u \ 
{\rm   to\ each\ line \  } L \ {\rm  is\ in\ }
 \overline{\rm C}{\rm onvex}(L\cap X)
\eqno{(\BB.6)}
$$

\Prop{\BB.6}  {\sl  Suppose $u\in\USC(X)$.  Then:
\smallskip
 \centerline{ $u\in \overline{\rm C}{\rm onvex}(X)$
  if and only if $u+v\in\SA(X)$ for all $v\in \SA(X)$.}
}

\pf
If $u\in {\rm Convex}(X)$ and $v\in \SA(X)$, then by Proposition \BB.5, $u+v\in\SA(X)$.
Furthermore, the extra case $u\equiv -\infty $ is obvious.

It remains to show that if $u+v\in\SA(X)$ for all $v\in\SA(X)$, then 
$u\in \overline{\rm C}{\rm onvex}(X)$.
It will suffice to show that: 
$$
u\notin \overline{\rm C}{\rm onvex}(X) \ \Rightarrow\ \exists {\rm \  a\
subaffine\ quadratic \ function\ } B {\ \rm with}\ u+B\notin \SA(X).
\eqno{(\BB.7)}
$$
Since $u\notin \overline{\rm C}{\rm onvex}(X)$, we know that for the restriction $\overline u$ 
of $u$ to some line $L$, we have $\overline u\notin \overline{\rm C}{\rm onvex}(L)$.
For $n=1$, $\SA  =\overline{\rm C}{\rm onvex}$, so that (\BB.3) applies to $\overline u$. 
Assume that the line  is the $x_1$-axis   and that the point  on the line $L$ in (\BB.3) is the origin
in $\rn$. Also, change $\overline u$ by the affine function in (\BB.3).  
Then there exists $\d>0$ so that
 $\overline u(t) \leq -\e t^2$ for $|t|\leq \d$ and $\overline u(0)=0$.
Hence, by the upper semicontinuity of $u$, there exists $r>0$ small with 
$$
\qquad\quad   u(t,y)+{\e\over 2} t^2 <0 \quad{\rm for\ \ } t=\pm\d, \ \  |y|\leq r.
$$
Now choose $\l >>0$ so that
$$
u(t,y) +\smfrac \e 2 t^2 - \l |y|^2\ <\ 0\quad{\rm for\ \ } |t|\leq \d,\ \  |y|=r.
$$
The quadratic function $B\equiv {\e\over 2}t^2 -\l |y|^2$ is subaffine by Proposition \BB.4, but the sum
$u+v$ is zero at the origin and strictly less than zero on the boundary of the cylinder $|t|\leq \d$,
$|y|\leq r$ about the origin.  Hence, $u+B$ is not subaffine.\qed


\vfill\eject

\noindent{\headfont \CC. Dirichlet Sets -- The Maximum of Two Functions.}  \smallskip

Each subset $F$ of $\Symn$ defines a class of $C^2$-functions $u$ by requiring
that $\Hess_x u\in F$ at each point $x$.  An important property that 
we want functions of this type $F$ to have is:
\medskip\noindent
{\bf The Maximum Property:}   If $u,v$ are of type $F$, then $\max\{u,v\}$ is of type $F$.
\medskip\noindent
Of course, we must first extend the definition of type $F$ functions to include non \ct-functions such as
$\max\{u,v\}$.  The appropriate condition on $F$ which insures this maximum property
is  the standard positivity (or elliptic) condition given in the next definition. 
See  Remark \CC.3 at the end of this section for more detail.

\Def{\CC.1}  A proper non-empty closed subset $F\ss \Symn$ will be called a {\sl Dirichlet set} if
it satisfies the {\sl Positivity Condition:}
$$
F+\cp\ \ss\ F
\eqno{(\CC.1)}
$$
where
$$
\cp \ \equiv \ \{A\in\Symn : A\geq0\}
\eqno{(\CC.2)}
$$
denotes the set of non-negative quadratic forms on $\rn$.
\medskip

We first introduce the notion of $F$-plurisubharmonicity for $C^2$-functions.
The definition will be substantially generalized in the next section.

\Def{\CC.2}  Suppose $F\ss\Symn$ is a Dirichlet set.  
 If $u\in C^2(X)$ has $\Hess_x u\in F$ for all $x\in X$, then $u$ is 
of   {\sl type $F$} or {\sl $F$-plurisubharmonic}.
 If   $\Hess_x u\in \Int F$ for some $x\in X$, then $u$ is called {\sl strict  of type $F$ at $x$}. 

\medskip
\noindent
{\bf Elementary Properties of Dirichlet Sets $F$:}
\smallskip
\item{(1)}\ \ $F+\Int\cp\ \ss\ \Int F$

\smallskip
\item{(2)}\ \ $F  \ =\   \overline{\Int F}$

\smallskip
\item{(3)}\ \ $\Int F+ \cp\ \ss\ \Int F$

\smallskip
\item{(4)}\ \ For each $B\in \Symn$ the set $\{t\in\bbr: B+tI\in F\}=[b,\infty)$ for some $b\in\bbr$.

\smallskip
\item{(5)}\ \ ($F$ is ``Asymptotically convex'') Given $A, B\in F$, 
$\exists \, t>0$ such that $A+tI$ and $B+tI$ 

\ \ both belong to the convex subset
$(A+\cp)\cap (B+\cp)$ of $F$.

\smallskip
\item{(6)}\ \ $F$ is Dirichlet $\Rightarrow$ $\l F+A$   is Dirichlet for $\l>0$ and $A\in \Symn$.

\smallskip
\item{(7)}\ \ $F$ is Dirichlet $\iff$ $gF$ is Dirichlet with $g\in GL_n(\bbr)$ acting on $\Symn$
by the

 \ \ standard action $g(A)=g^t \circ A\circ g$.

\medskip
\noindent
{\bf Proofs:}
\smallskip
\item{(1)} For each $A\in F$ the open set $A+\Int \cp$ is contained in $F$.

\smallskip
\item{(2)} Use (1) and  $A=\lim_{\e\to0}(A+\e I)$.

\smallskip
\item{(3)} For each $P\in \cp$ the open set $\Int F+P$ is contained in $F$.

\smallskip
\item{(4)} Since $F$ Dirichlet implies that $F-B$ is Dirichlet, we may assume that $B=0$.
We must show that the set $\L_F \equiv \{t\in \bbr :   tI\in F\}$ is connected, proper,
 and non-empty.  If $t_0\in \L_F$, then by the Positivity Condition $t\geq t_0$ implies $t\in \L_F$. Hence,  $\L_F$ is connected.
 If $\L_F=\bbr$, then $-tI\in F$ for all $t>0$.  Hence, $-tI+\cp\ss F$  for all $t>0$.
 This implies that $F$   equals $\Symn$ which is not allowed. Therefore, $\L_F\neq \bbr$.
 This implies $\L_F\neq \emptyset$ by duality. (See Remark \DD.2 in the next section.)

\smallskip
\item{(5)} Pick $t>>0$ so large   that $A+tI\in B+\cp$ and $B+tI\in A+\cp$.

\smallskip
\item{(6)} and (7) are straightforward.

\Remark{\CC.3}  Motivation for the Positivity Condition  is provided by
\medskip\noindent
{\bf The Hessian Lemma:}  {\sl Suppose  $u,v\in C^2(X)$ and $\n(u-v)\neq 0$ on $\{u=v\}$. 
Then taking the distributional hessian, we have
$$
\Hess(\max\{u,v\}) \ =\ \chi_{\{u\geq v\}} \Hess u  + \chi_{\{v\geq u\}} \Hess v + \mu \n(u-v)\circ
\n(u-v)
$$
where $\mu$ is a non-negative measure supported on $\{u=v\}$.}
\medskip\noindent
This formula strongly suggests that one should require:
$$
A+\x\circ\x \in F \qquad {\rm for\ all\ }  A\in F, \ \x\in\rn.
$$
Since each $P\geq 0$ can be written as $P=\sum_j \l_j e_j\circ e_j$,  this condition
is equivalent to the Positivity Condition (\CC.1) that $F+\cp\ss F$.
We omit the proof of this lemma.
\medskip

\vfill\eject


\noindent{\headfont \DD. Dirichlet Duality.}  \smallskip

As noted in Definition \CC.2, each subset $F$ of $\Symn$ defines a class of $C^2$-functions 
$u$ by requiring that $\Hess_x u\in F$ at each point $x$.
 In this section we will give a dual characterization
of this condition, which will enable us to define functions of type $F$ which are not necessarily 
of class \ct.  This nonlinear duality can be used in a fashion which has some similarity  to the use of 
distribution theory in linear problems.

Throughout  this section we assume that $F$ is a Dirichlet set.
Let
$$
\cpt \ =\ \sim (-\Int \cp) \ =\ -(\sim\Int \cp).
$$
denote the set of all quadratic forms except those that are negative definite, i.e., $A\in\cpt$ iff $A$ has at least one eigenvlaue $\geq0$.  Note that for $u\in C^2(X)$,
\medskip

\centerline{ $u$ is convex iff $u$ is of type $\cp$
\qquad and\qquad 
$u$ is  is subaffine iff $u$ is of type $\cpt$.}
\medskip

The second statement is just Proposition \BB.4

The key to the dual characterization of functions of type $F$ is the existence of a dual subset
$\wt F$.  This is made precise in Lemma \DD.3  below.

\Def{\DD.1}   Suppose $F\ss\Symn$ is a Dirichlet set.  The {\sl Dirichlet dual} of $F$
is the set
$$
\ft \ =\ \sim (-\Int F) \ =\ -(\sim\Int F).
$$
{\bf Elementary Properties of the Dirichlet Dual.}
\smallskip
\item{(1)}\ \ $\wt{\wt{F}}\ =\ F$. \hskip 1.8in     (3)\ \ $\wt{F_1\cap F_2}\ =\ \ft_1\cup \ft_2$

\smallskip
\item{(2)}\ \  $F_1\ \ss\ F_2\quad\Rightarrow\quad \ft_2\ss\ft_1$.
\qquad\qquad\quad\  (4)\ \ $\wt{F_1\cup F_2}\ =\ \ft_1\cap \ft_2$

\smallskip
\item{(5)}\ \ $\wt{F+A}\ =\ \ft - A$.

\smallskip
\item{(6)}\ \ $F$ is a Dirichlet set \ $\iff$\  $\ft$ is a Dirichlet set.

\medskip
\noindent{\bf Proofs.}
\smallskip
\item{(1)} follows from $F=\overline{\Int F}$.

\smallskip
\item{(2)} (3) and (4) are obvious.

\smallskip
\item{(5)} Note that $B\in \wt{F+A} \iff -B \notin \Int(F+A) = \Int F +A
\iff -(B+A)\notin \Int F \iff B+A\in \ft \iff B\in\ft-A$.

\smallskip
\item{(6)} Suppose $P\in \cp$.  Then   $F+P\ss F$ or equivalently $F\ss F-P$.  By (2) this implies that 
$\wt{F-P}\ss\ft$.  By (5) we have $\wt{F-P} = \ft+P$ so that $\ft+P\ss\ft$.

\Remark{\DD.2} Define $\wt{\L}_F = \sim(-\Int \L_F)$ and note that $\wt{\L}_F = \L_{\ft}$.
Hence $\L_F =\emptyset \Rightarrow \L_{\ft}=\bbr \Rightarrow \ft=\Symn \Rightarrow F=\emptyset$,
completing the proof of Property (4) in Section \CC.
\medskip

The following duality result is stated in several  forms:  first for the special case of points  $A\in\Symn$ (i.e., quadratic functions), and then for functions  $u\in C^2(X)$.

\eject

\Lemma{\DD.3}  {\sl  Suppose $F$ is a Dirichlet set.  Then
\medskip

(1)\ \ \qquad\qquad\qquad\quad\ 
$A\in F \quad \Leftrightarrow\quad A+B\in\cpt \quad$ for all $B\in \ft$.

\medskip

(2) \ \ $u\in C^2(X)$ is of type $F$ $\ \  \Leftrightarrow\quad u+B\in\SA(X)$\ \ 
  for all quadratic  $B\in\ft$.

\medskip

(3) \ \ $u\in C^2(X)$ is of type $F$ $\quad \Leftrightarrow\quad u+v\in\SA(X)$\ \ 
  for all $v\in C^2(X)$ of type $\ft$.
  }

\pf  Statement (3) follows from the special case (1) by setting $A=\Hess_x u$,
$ B=\Hess_x v$, and using Definition \CC.2 along with Proposition \BB.4. 
Thus the three conditions are equivalent.

To prove (1), first note that:\smallskip
\item{(1)$'$} \ \ $A\in F \iff A+\cp\ss F$
\smallskip
\noindent
is obviously true because of the positivity condition (\CC.1).

Now $A+\cp \ss F \iff \ft\ss\wt{    A+\cp}$ (which equals $\cpt-A$) $\iff A+\ft\ss\cpt$.  Thus (1)$'$ is equivalent to:\smallskip
\item{(1)} \ \ $A\in F \iff A+\ft\ss \cpt$.
\qed

\medskip

Because of this Lemma we can extend our Definition \CC.2 of type $F$ 
from \ct-functions to upper semi-continuous functions.This extension 
is another central concept of the paper.

\Def{\DD.4}  A function $u\in \USC(X)$ is said to be of {\sl type $F$} or {\sl $F$-\psh} if 
$$
u+v \in \SA(X) \fa v\in C^2(X) \ \ {\rm of \ type\ }\ft.
\eqno{(\DD.1)}
$$
Let $F(X)$ denote the set of all $u\in \USC(X)$ of type $F$.

\Prop{\DD.5}  {\sl  Suppose $u\in \USC(X)$.  Then (for $X$ connected)
$$
\eqalign
{
u \ {\sl is\ convex\ or\ } u\equiv-\infty\ \ \  &\iff \ \ \ u\ \  {\sl is\ of\ type\ }\cp
\qquad{\sl and}   \cr
u \ {\sl is\ subaffine\ } \ \ &\iff \ \ \ u \ \ {\sl is\ of\ type\ }\cpt.
}
$$
Moreover,  for any $u$ of type $\cp$ and any $v$  of type $\cpt$, the sum  $u+v\in\SA(X)$.
}
\pf
 This is just a restatement of Propositions \BB.5 and \BB.6.\qed\medskip
 
 Note that for two Dirichlet sets $F_1$ and $F_2$,
 $$
 F_1(X)\ \ss\ F_2(X) \qquad\iff \qquad F_1\ \ss\ F_2
 \eqno{(\DD.2)}
 $$

It is important to have some equivalent formulations of the definition of
functions of type $F$.   For example, as it stands it is not clear that if $u$ is
of type $F$ on $X$, then the restriction of $u$ to a smaller open subset is 
also of type $F$. This however is true and is easily seen from other equivalent 
definitions.

In making these reformulations we first reduce the space of test functions
from $C^2(X)\cap \ft(X)$ to just $\ft$, the space of quadratic functions of type $F$.
The second formulation says that if $u\notin F(X)$, then near some point 
$x_0\in X$, the condition for type $F$ is strongly violated.

\Lemma{\DD.6}  {\sl  A function $u$ is in $F(X)$ if and only if 
$$
u+B\in\SA(X) \quad   {\sl for\ all\ quadratic\ functions\ } B\in \ft.
\eqno{(\DD.3)}
$$
Moreover, $u\notin F(X)$ if and only if 
$$
\eqalign{
& {\sl there\ exist\ } B\in \Int\ft, x_0\in X, a\ {\sl affine,\ and\ } \e>0
\ {\sl such\ that\ } 
 \cr
&u+B-a\ \leq \ -\e|x-x_0|^2 \ \ {\sl near \ }x_0\ \ {\sl and }\cr
&\qquad\quad\quad\ \  =0 \qquad\qquad \qquad{\rm  at\ } x_0.
}
\eqno{(\DD.4)}
$$
}

\pf
If $u\in F(X)$, then, taking $v=B$, we see that (\DD.1) implies (\DD.3).
Furthermore, (\DD.3) obviously implies that (\DD.4) is false.
It remains to show that if (\DD.1) is false, then (\DD.4) is true.
If (\DD.1) is false, then there exists $v\in C^2(X)\cap \ft(X)$ such that 
$u+v\notin\SA(X)$.  Applying Lemma \BB.2, there exist
$x_0\in X, \e>0$ and and affine function $a$ with 
$u+v-a\leq -2\e|x-x_0|^2$ near $x_0$ and equal to zero at $x_0$.
Since $v\in C^2(X)$,replacing $v$ by the quadratic part $B$ of $v$ at $x_0$ yields:
$u+B-a\leq -2\e|x-x_0|^2$ near $x_0$ and $u+B-a=0$ at $x_0$.
Finally, since $B\in\ft$, we have $ B+\e I\in \Int\ft$, proving (\DD.4).
\qed

\bigskip
\centerline
{\bf Properties of the class  $F(X)$ for Dirichlet Sets $F$.}

\medskip

\item {(1)}  (Local Property). A function $u$ is locally of type $F$ if and only if  $u$ is  (globally)   of type $F$.

\medskip

\item {(2)} (Affine Property). $F(X)+{\rm Aff}(X)\ss F(X)$, i.e., if $u\in F(X)$ and $a$ is affine, then $u+a\in F(X)$.

\medskip

\item {(3)} (Translation Property). If $u\in F(X)$, then $v(x)\equiv u(x-y) \in F(X+y)$.

\medskip

As anticipated, the Positivity Condition insures the maximum property.

\medskip

\item {(4)} (Maximum Property). If  $u, v \in F(X)$,  then $\max\{u,v\}\in F(X)$.

\medskip

\item {(5)}   (Decreasing Limits) If $\{u_j\}_{j=0}^\infty$ is a decreasing (i.e., $u_j\geq u_{j+1}$)  sequence of functions  in $F (X)$, then  $\lim_j u_j\in F (X)$.

 \medskip

\item {(5)$'$}   (Uniform Limits) If $\{u_j\}_{j=0}^\infty$ is a  sequence of functions  in $F (X)$
which converges uniformly to $u$ on compact subsets, then  $u\in F (X)$.

 \medskip

\item {(6)}   (Families Locally Bounded Above) Suppose $\cf\ss F(X)$ is 
locally uniformly bounded above.  Then the upper envelope $u=\sup_{f\in\cf} f$ has u.s.c. regularization
$u^*\in F(X)$.

 \medskip

\item {(7)}   If $u$ is twice differentiable at $x\in X$, then $\Hess_x u\in F$.

\medskip\noindent
{\bf Proofs.}\smallskip
\item{(1)} By Definition \BB.2, subaffine functions restrict to subaffine functions, and by Lemma \BB.2, 
locally subaffine implies subaffine.
Condition   (\DD.3) now ensures that functions of type $F$ are locally of type $F$
(since the quadratic functions are ``universal'', i.e., defined on all of $\rn$).
Conversely, if $u$ is locally of type 
$F$, then (\DD.4) is false, and hence $u$ is globally of type $F$.
\smallskip
\item{(2)}, (3) and (4) follow from the definitions.
 \smallskip
 \item{(5)} This standard proof uses the fact that 
 $\{x\in \partial K: u_j(x)+v(x)\geq a(x)+\e\}$ is compact for $u_j\in\USC(X)$ and hence empty for $j$ large.
  \smallskip
 \item{(5)$'$}  This is standard from (5).  Given $\e_j\searrow 0$, $\e_{j+1} <\e_j$, pick $j$ large so that 
 $|u_j-u| < \half(\e_{j+1}-\e_j)$ and set $u_j' = u_j+ \half(\e_{j+1}+\e_j)$.
 
 \smallskip
 \item{(6)} 
 If $u^*+v\leq a$ on $\partial K$, then $f+v\leq a$ on $\partial K$ for all $f\in \cf$.  Hence,  $f+v\leq a$ on $K$ for all $f\in \cf$, and so $u+v\leq a$ on $K$. Since $v\in C^2$, $u^*+v=(u+v)^*\leq a^*=a$
 on $K$.
  \smallskip
 \item{(7)} 
 Let $H=\Hess_x u$.  Then the quadratic function $H(y)$ is the uniform limit 
 as $\e\to0$ of the approximate hessians
 $$
 H_\e(y)\ =\ \e^{-2}\{u(x+\e y)-u(x)- \n u(x)\cdot y\}.
 $$
 By (5)$'$ it suffices to prove that:
 $$
 {\rm The \  approximate\ hessians \ }  H_\e\ {\rm are\ of\ type\ } F
 \eqno{(\DD.5)}
 $$
 Since ${1\over \e^2}(u(x)+\e  \n u(x)\cdot y )$ is affine,
 we must show that $(Lu)(y) = {1\over \e^2}u(x+\e y)$ is of type $F$. 
 Now $L$ has inverse given by
 $$
( L^{-1} v)(y) \ =\ \e^2 v\left({y-x\over \e}\right).
 $$
Note that if $v$ is \ct, then $\Hess\, L^{-1} v = \Hess\, v$ at corresponding points.
Consequently,   if $v$  is \ct  \ and of type $\ft$, then $L^{-1} v$ is of type $\ft$.
Therefore,  $u+L^{-1} v  \in\SA(X)$ because $u\in F(X)$.
  Finally, since $L$ maps subaffine to subaffine 
  (This is Property (7) for $\cpt$ and can be verified directly), 
  we conclude that $Lu+v = L(u+L^{-1} v)$ is subaffine.     
Hence, $Lu$ is of type $F$ as desired.

\medskip

\Remark{\DD.7. (The Maximum Principle)} This principle is not always true and (perhaps 
surprisingly) not necessary for the Dirichlet problem.  However, the Maximum Principle
for all functions in $F(X)$ is true if and only if $F(X)\ss\SA(X)$ because of (2) above.  
By (\DD.2) this is equivalent to $F\ss \cpt$.  Note that
$$
F \ss \cpt\qquad\iff \qquad 0\notin \Int F
\eqno{(\DD.6)}
$$
If $0\in\Int F$, then $F$ contains a negative definite quadratic form so that
$F\ss\cpt$ is impossible. 
Conversely, if $F\not\ss \cpt$, then $F$ contains a negative definite quadratic form $A=-P_0$,
$P_0>0$.  The open set $\{A+P:P>0\}\ss F$
contains the origin.  This proves (\DD.6) and hence
\Prop{\DD.8} {\sl Suppose $F\ss\Symn$ is a Dirichlet set.  The maximum principle holds for each 
$u\in F(X)$  if and only if $0\notin \Int F$.}
\medskip

In the cases where $0\in \Int F$, (i.e., when the maximum principle does not hold), the functions $u\in F(X)$ will be called $F$-{\sl quasi-plurisubharmonic}.

\Remark{\DD.9. (Viscosity Subsolutions)}  The condition (\DD.4) above is equivalent to:
$$
\eqalign{
 \exists\, x_0\in &X   \ \ {\rm and}\ \ \psi\in C^2(X) \ \ {\rm which\ is\ strict\ of\ type\ } \ft \ {\rm at}\  x_0
 \cr
\ &{\sl such\ that\ } \ u+\psi \ \  {\rm  has\ a\ local\ maximum\ at\ } x_0.
}
\eqno{(\DD.4)'}
$$
{\bf  Proof.}  That  (\DD.4) $\Rightarrow$ (\DD.4)$'$ is obvious.
For the converse set $-a = \langle(\n \psi)(x_0), x-x_0\rangle$ and $B=\half
(\Hess_{x_0}\psi)(x-x_0)-2\e I\in \Int \ft$. \qed
\smallskip

Since $\Int\ft = -(\sim F)$, if we set $\vf = -\psi$, then (\DD.4)$'$ is equivalent to 
the condition:
$$
\eqalign{
 \exists\, x_0\in X   \ \ {\rm and}\ \ \vf \in C^2(X) \ \ &{\rm with\ }\ \Hess_{x_0}\vf  \notin F
 \cr
\ {\rm but \ } \ u-\vf \ \  {\rm  has\ a\ local}\ &{\rm maximum\ at\ } x_0.
\cr
}
\eqno{(\sim V)}
$$
Finally the negation of ($\sim$ V) is:
$$
\eqalign{
&\qquad\forall \, x_0\in  X   \ \ {\rm and}\ \ \vf \in C^2(X) \ \  
 \cr
\ {\sl if \ } \ u-\vf \ \  &{\rm  has\ a\ local\ maximum\ at\ } x_0, \ {\rm then\ } \Hess_{x_0}\vf  \in F
}
\eqno{(V)}
$$
Condition (V) is the standard viscosity definition of subsolution.


\vfill\eject

\noindent{\headfont \EE. Boundary Convexity.} \medskip

We assume throughout this section that $\O$
 is a bounded domain in $\rn$  with  smooth boundary $\bo$.
It turns out that  the natural boundary-convexity condition  associated
to a Dirichlet set  $F$ is expressed in terms of another \ds $\oa F$, 
the {\sl ray set  associated with $F$}, which will be defined
in a moment.  
Our main result (Theorem \EE.12 below)
asserts that  strict local $\Fa$-convexity of the boundary
implies the existence of a global  defining function which is  strictly $\Fa$-plurisubharmonic.
 This global function  will play a key role in our solution to the Dirichlet problem in \S \FF.

The key property of this associated ray set $\Fa$ is the

\medskip
\noindent
{\bf Ray Property:}
$$
A\in \Fa \quad \iff \quad tA\in \Fa \fa  t\geq 0.
\eqno{(\EE.1)}
$$
Moreover, $\Fa = F$ if and only if $F$ itself has the ray property.
\medskip

A \dir set $\Fa$ with property (\EE.1) will be called a {\bf Dirichlet-Ray set} or {\bf D-Ray set}.
We assume for the moment that $\Fa$ is any Dirichlet-Ray set (not necessarily the one 
associated with $F$).
Note that since scalar multiplication by $t>0$ is a homeomorphism of $\Symn$, (\EE.1) implies  
$$
A\in\Int\Fa \quad\iff\quad tA\in \Int \Fa\fa t>0.
\eqno{(\EE.1)'}
$$

A smooth function $\rho$ defined near a point $x\in\bo$ is said to be a 
{\sl local defining function for $\bo$ near $x$} if  on some neighborhood of 
$x$, we have $\O=\{\rho<0\}$ and $\n \rho\neq 0$. At the boundary point
$x$, let $T=T_x\bo$ denote the tangent space and $n$ a unit normal vector.

\Def{\EE.1} The boundary $\bo$ is  {\sl    strictly $\Fa$-convex at a point $x\in \bo$}
if
$$
\Hess_x \rho\bigr|_T\ =\ B\bigr|_T \ \ {\rm  for\  some\ } B\in \Int \Fa
\eqno{(\EE.2)}
$$

 \Lemma {\EE.2} {\sl The  condition of strict $\Fa$-convexity for $\bo$ is 
 independent of the defining function $\rho$.}
 
 \pf Any other defining function $\wt \rho$ is of the form $\wt \rho = u\rho$ with $u>0$.
 At $x\in \bo$, $\Hess_x \wt\rho = u \Hess_x \rho + \n u\circ\n \rho$.  Since $\n u\circ \n \rho $ restricts to be zero on $T$, we have $\wt H\bigr|_T = u H\bigr|_T$. By the ray property (\EE.1)$'$ for $\Int \Fa$ the proof is complete.\qed \medskip

The notion of strict $\Fa$-convexity has other   useful formulations.
Let $P_n= n\circ n \in\Symn$ denote orthogonal projection onto the line in the normal direction $n$

 \Lemma{\EE.3}  {\sl The following conditions on a local defining function
 $\rho$ for $\bo$  are equivalent.
 \medskip
  
  (1) \ \ $\Hess_x \rho\bigr|_T\ =\ B\bigr|_T$ for some $B\in \Int \Fa$ (i.e. $\bo$ is  strictly $\Fa$-convex at $x$).
 
 \medskip
  
  (2) \ \ $\Hess_x \rho\bigr|_T +t P_n\ \in\  \Int \Fa$ for all $t \geq $ some $t_0$.
 
 \medskip
  
  (2)$'$ \ \ $\Hess_x \rho  +t P_n\ \in\  \Int \Fa$ for all $t \geq $ some $t_0$.
  }
  
\pf
Let $H=\Hess_x\rho$.
Statements (2) and (2)$'$  each imply (1) since in both cases the restriction
to $T$ equals $H\bigr|_T$.
   Suppose now that (1) is true.  Then, in terms of the 2$\times$2 blocking induced by
 $\rn=\span n \oplus T$, we have 
 $
H-B\ =\ \left( \matrix{a & \a\cr\a&0\cr}\right).
 $
 Therefore, 
 $
 H+tP_n \ =\ B-\e I+tP_n+H-B+\e I.
 $
 If $\e>0$ is chosen small, then $B-\e I \in \Int \Fa$, while
 $
 tP_n +H-B+\e I \ =\ \left(\matrix{t+a+\e & \a\cr \a &\e I\cr} \right)
 $
 which is positive definite for $t>>0$.
 By Property  (3) for \dir sets, this implies that $H+tP_n \in \Int \Fa$ for $t>>0$. 
 \qed
 
\Cor{\EE.4} {\sl  Let $II$ denote the second fundamental form of $\bo$ with respect to the 
inward-pointing unit normal $n$. Then $\bo$ is strictly $\Fa$-convex at $x\in\bo$ if and only if 
 \medskip
  
  (1) \ \ $II_x  \ =\ B\bigr|_T$ for some $B\in \Int \Fa$,
  or equivalently
 
 \medskip
  
  (2) \ \ $II_x  +t P_n\ \in\  \Int \Fa$ for all $t \geq $ some $t_0$.
 
}

\pf
By Lemma \EE.2 we  may choose $\rho$ to be the signed distance function in a neighborhood
of $\bo$,  i.e.,   for $x$ near $\bo$
$$
\rho(x)\ =\ \d(x)\ =\ \cases{-\dist(x,\bo)\quad {\rm if\ \ } x\in \O  \cr
+\dist(x,\bo)\quad {\rm if\ \ } x\notin \O  \cr}.
$$
 Then it is a standard calculation (cf. [HL$_2$, (4.7)]) that 
 $$
 \Hess\, \d \ =\ \left(  \matrix{
 0&0\cr
 0& II
 }
 \right).
 $$
 with respect to the splitting $\rn=(\bbr \cdot \n \d) \oplus (\n \d)^\perp$.
We now apply Lemma \EE.3.\qed

\bigskip
\noindent
{\bf The Associated Ray Set  $\oa F$ of $F$.}

\Def{\EE.5}
Suppose $F$ is a Dirichlet set and $B\in \Symn$ is fixed.  The {\sl ray set with vertex
$B$ associated to $F$}, denoted by $\oa {{F_B}}$, is defined by
$$
\oa {{F_B}}\ =\ \{A\in\Symn : \ \exists \,  t_0\ \ {\rm such \ that\ }\  B+tA\in F  \fa t\geq t_0\}
$$

\Ex{\EE.6}  The set $\oa {{F_B}}$ may not be closed.  Take $B=0$ and $F=\cp\cap\{\det\geq1\}$.

\Lemma{\EE.7} {\sl The closure of $\oa {{F_B}}$ is independent of the vertex $B$.}

\pf  This is Property (6) proven below.

\Def{\EE.8}  Suppose $F$ is a Dirichlet set. The {\sl ray set associated to $F$}, denoted
$\oa F$, is defined to be the closure of $\oa {{F_B}}$ for any vertex $B$.

\def\Closure{{\rm Closure}}

\medskip
\noindent
{\bf Elementary Properties of $\oa {{F_B}}$:}
\smallskip
\item{(1)} \ \ $\oa {{F_B}}+\cp\ \ss\ \oa {{F_B}}$

\smallskip
\item{(2)} \ \ $\oa {{F_B}}+\Int\cp\ \ss\ \Int \oa {{F_B}}$

\smallskip
\item{(3)} \ \ $\oa {{F_B}}\ \ss\ \Closure( {\Int \oa {{F_B}}})$

\smallskip
\item{(4)} \ \  $\Int \oa {{F_B}} \ \ss\ \oa {{F_{B'}}}$\ \  for all $B'$

\smallskip
\item{(5)} \ \ $\Int \oa {{F_B}} \ =\ \Int \oa {{F_{B'}}}$\ \  for all $B'$

\smallskip
\item{(6)} \ \ $\Closure( \oa {{F_B}}) \ =\ \Closure( \oa {{F_{B'}}})$ \ \ for all $B'$

\medskip
\noindent
{\bf Proofs.}
\smallskip
\item{(1)} $B+tA\in F$ for $t\geq t_0$\  $\Rightarrow\ B+t(A+P)=B+tA+tP\in F$ for $t\geq t_0\geq0$.

\smallskip
\item{(2)} If $A\in \oa {{F_B}}$, then by (1) the open set $A+\Int\cp \ss\oa {{F_B}}$, and hence
$A+\Int\cp \ss\Int\oa {{F_B}}$.

\smallskip
\item{(3)} If $A\in  \oa {{F_B}}$, then by (2) we have that for $\e>0$, 
$A+\e I\in\Int  \oa {{F_B}}$. Hence, $A=\lim_{\e\to0}(A+\e I) \in \overline{\Int \oa {{F_B}}}$.
Note that for Example 1 equality in (3) does not hold.

\smallskip
\item{(4)}   If $A\in \Int  \oa {{F_B}}$, then  $A-\e I \in \oa {{F_B}}$ for some $\e>0$.
This means that there exists a $t_0$ so that $B+t(A-\e I) \in F$ for all $t\geq t_0$.
Now $B'+tA = B+t(A-\e I) + t\e I -(B-B')$.  Choose $\l>0$ so that $\l I -(B-B') >0$ is positive
definite.  If $t\geq t_0$ and $t\geq {\l\over \e}$, then $B'+tA \in F+\cp\ss F$, proving that
$A\in \oa {{F_{B'}}}$.

\smallskip
\item{(5)}  follows from (4)

\smallskip
\item{(6)}  follows from  (5) and (3).

 \qed
 \medskip
 
 Since boundary convexity involves $\Int \Fa$, some additional facts about $\Int \Fa$
 are useful.  The associated ray set $\Fa$ for $F$ was defined to be as large as possible.
 The smallest set of rays associated with $F$ is $\Int \oa{{(\Int F )_B}}$ where
 $$
 \oa{{(\Int F )_B}}\ \equiv \ \{A\in\Symn : \exists\  t_0\ {\rm so\  that\ }
 B+tA\in \Int F\ \  \forall \ t\geq t_0\}.
 $$ 
 \Ex{} The set $ \oa{{(\Int F )_B}}$  may not be open.  Take $F=\cp$ and $B=I$.

 \Lemma{\EE.9}
 $$
 \Int  \oa{{(\Int F )_B}}\ =\ \Int \Fa.
 $$
 \pf
Since $\oa{{(\Int F )_B}} =  \oa{{(\Int (F -B))_0}}$ and $\oa{{F-B}} = \Fa-B$, we may assume $B=0$
Because $\oa{{(\Int F )_0}}\ss \Fa$,
  it suffices to show 
 $\Int \Fa \ss \oa{{(\Int F )_0}}$.
 Suppose $A\in \Int\Fa$.  Then there exists $\e>0$ with $A-\e I\in \Fa = {\rm Closure}(\oa{{ F_0}})$.
 Therefore for all $\d>0$ there exists $B\in  \oa {{F_0}}$ with
 $|A-\e I-B  | <\d$, which implies that $\d I +A -\e I -B >0$.
 Take $\d={\e\over 2}$. Then $A -{\e\over 2} I -B >0$ and $B\in  \oa {{F_0}}$.
 Hence, there exists $t_0>0$ so that $t\geq t_0\ \Rightarrow\ tB\in F$.
 Therefore, $t(B+{\e\over 2} I ) = tB + {\e t\over 2} I \in \Int F$ if $t\geq t_0$.
 This proves that $B+{\e\over 2} I  \in \oa{{(\Int F )_0}}$.
 Finally, $A=B+{\e\over 2} I  +A -B-{\e\over 2} I  \in \oa{{(\Int F )_0}} +\Int \cp\ss 
 \oa{{(\Int F )_0}}$.  \qed

\Cor{\EE.10}  {\sl One has $A\in\Int\Fa$ if and only if }
$$
 \exists\  \e>0\ {\sl and \ } R>0\ {\sl such\ that\ } C(A-\e I)\in F\ \ {\sl for\ all\ }C\geq R
\eqno{(\EE.3)}
$$ 
\pf
It follows easily from the previous lemma that $\Int \Fa=\Int \oa{{F_0}}$.
Therefore, 
 $A\in\Int\Fa$, then for some $\e>0$, $A-\e I\in \oa{{F_0}}$, i.e., there exists $R>0$ such that $C\geq R$ implies $C(A-\e I)\in F$.

Conversely, assume (\EE.3) is true.  Suppose $B+\e I>0$. This condition defines a neighborhood of the origin in $\Symn$.  It suffices to show that $A+B\in \Fa$ for all such $B$.  Now $A+B =A-\e I +B +\e I$ and hence
$C(A+B) = C(A-\e I) + C(B+\e I)$ which belongs to $F+\cp\ss F$ if $C\geq R$.  Hence
$A+B\in \oa{{F_0}}\ss \Fa$.\qed

 \Prop{\EE.11}  {\sl If $F$ is a Dirichlet set, then the associated ray set $\oa F$ is also a 
 Dirichlet set.   Moreover,  $\Fa$ has the ray property.}
 
 \pf Take the closure in (1) above.\qed

\Remark{}  Since the associated ray set $\Fa$ of $F$ is a Dirichlet-Ray set,
 the definition of strict $\Fa$-convexity at a boundary point $x\in \bo$ is independent of the defining function $\rho$ (Lemma \EE.2).
\medskip

A smooth function $\rho\in C^\infty(\Ob)$ is called a  {\sl global defining function for $\bo$}
if   $\O=\{\rho<0\}$ and $\n\rho\neq0$ on
$\bo$.

\Theorem{\EE.12}  {\sl Suppose $\Fa$ is a Dirichlet-Ray set.
 If the boundary $\bo$ is strictly $\Fa$-convex at each point,
then there exists  a global defining function $\rho \in C^\infty(\Ob)$
 for $\bo$  which is strict of type $\Fa$
on $\Ob$.  Moreover, if $\Fa$ is the ray set associated with a \dir set $F$, then}
$$
 \exists\ \e>0 \  {\rm and }  \ R>0 \ {\rm such\  that\ }\  C\left(  \rho - \e \half |x|^2\right)\ \in \   F(\Ob) \fa C\geq R
 \eqno{(\EE.4)}
 $$

The existence of the function $\rho$ in this theorem is the only part of this section needed to
solve the \dir Problem in the \S \FF.

 \pf Pick any smooth defining function $\rho \in C^\infty(\Ob)$  for  $\bo$.  
 Let $\wt \rho = \rho +C\rho^2$. 
 Since $\bo$ is strictly $\Fa$-convex at each point,  we have by
  Lemma \EE.3,  part (2)$'$,  that on $\bo$, 
 $\Hess \wt \rho = (1+2C\rho) \Hess\rho +C\n \rho\circ\n\rho
  = \Hess\rho +C\n \rho\circ\n\rho \in \Int \Fa$ for all $C>>0$.
That is,  for large $C$, the defining function $\wt \rho$ is strictly $\Fa$-plurisubharmonic at each 
boundary point. 
This proves that we may assume  the defining function $\rho$ is strict of type $\Fa$ 
 in a neighborhood of $\bo$.  Choose $r >0$ so that the set
$ \{  -r <  \rho  < r\}$ is contained in this neighborhood
of $\bo$ where $\rho$ is strict.  Choose $\d>0$ small enough so that $-r+\d|x|^2 < \rho$
in a neighborhood $U$ of $\bo$.
We extend $ \rho$ to $\Ob$ by  setting
$$
\wh \rho \ \equiv\ \max\{ \rho, -r+ \d  |x|^2\}.
$$
On the open set $\O_{-r} = \{ \rho< -r\}$ we have $\wh \rho =  -r+ \d  |x|^2$, while on 
the neighborhood $U$ of $\bo$,  $\wh \rho = \rho$.
 Therefore,
 by the Maximum Property (4) of Section 4,  $\rho$ is strict of type $\Fa$ on $\Ob$.
   
   To complete the proof we smooth the maximum
    $\wh \rho \equiv M(u_1,u_2)\equiv \max\{u_1,u_2\}$
   of $u_1=\rho$ and $u_2= -r+\d |x|^2$, without changing $M(u_1,u_2)$ on the set where 
   $|u_1-u_2|\geq \e$.  Then choosing $\e>0$ small enough, the smoothing ${\wh\rho}_\e$
   will equal $\wh\rho$ in a neighborhood of $\bo$.
   Let  $M_\e(t_1,t_2)$ denote the smoothing of $M(t_1,t_2) = \max\{ t_1, t_2\}$ on $\bbr^2$
   (see [HL$_2$, Remark 1.6] for more details).  This can be done so that:
 
   \medskip
   \item{(1)} \ \ $M_\e(t_1,t_2) = M(t_1,t_2)$\ \ \  if $|t_1-t_2|\geq \e$.
  
   \medskip
   \item{(2)} \ \ ${\partial M_\e \over \partial t_1} +  {\partial M_\e \over \partial t_2} \ =\ 1$,
   ${\partial M_\e \over \partial t_1} \ \geq \ 0$,
  ${\partial M_\e \over \partial t_2} \ \geq \ 0$.
 
    \medskip
   \item{(3)} \ \ $M_\e(t_1,t_2)$ converges uniformly to $M(t_1,t_2)$ as $\e\to 0$.
   \medskip

It remains to show that ${\wh \rho}_\e   = M_\e(u_1,u_2)$ is strict of type $\Fa$  at each point
$x\in \Ob$.  By the chain rule
$$
\Hess M_\e(u_1,u_2)\ =\ {\partial M_\e \over \partial t_1}\,  \Hess\, u_1 + 
{\partial M_\e \over \partial t_2} \, \Hess\,  u_2 + 
\sum_{i,j=1}^2 {\partial^2 M_\e \over \partial t_i\partial t_j} \, \n u_i\circ \n u_j.
$$
One can show that the third term is $\geq 0$.  Hence,  by (2)  above it suffices to show  that
$$
A_s \ =\ s\Hess_x \rho + (1-s)2\d I\ \in\ \Int \Fa.
$$
At all points in the neighborhood of $\bo$ where $s={\partial M\over \partial t_1}\neq 0$, 
we have $\Hess_x\rho \in \Int \Fa$ and hence, $A_s\in \Int\Fa$.  At points $x$ where $s=0$, 
$\d I \in \Int \Fa$.
\qed

\bigskip

We conclude this section by listing some of the properties
 of Dirichlet-Ray sets and their corresponding  \psh functions.

\medskip\noindent
{\bf Elementary Properties of Dirichlet-Ray Sets $F$:}
\smallskip

\item{(1)}\ \    $\{t\in \bbr : tI\in F\}=[0,\infty).$

\smallskip

\item{(2)}\ \  $0\in \partial F$ \qquad\qquad\qquad\qquad (2)$'$ \ \ $0\in \partial\ft$

\smallskip

\item{(3)}\ \  $\cp\ss  F$ \qquad\qquad\qquad\qquad (3)$'$ \ \ $\ft \ss\cpt$

\smallskip

\item{(4)}\ \  $A\in \Int F$ \qquad $\iff$\qquad $tA\in \Int F \fa t\geq0$.

\smallskip

\item{(5)}\ \  $F$ is a D-ray set   \ $\iff$\ $\ft$ is a D-ray set

\smallskip

\item{(6)}\ \  $F\ \ss\cpt$ \qquad\qquad\qquad\qquad (6)$'$ \ \ $\cp\ \ss \ \ft$
\medskip

\noindent
{\bf Proofs:} Since $0\in F$ and $F+\cp\ss F$, we have $tI\in F$ for all $t\geq0$.
If $-tI\in F$ for some $t>0$, then $-tI+\cp\ss F$ for all $t>0$ and $F=\Symn$ contrary 
to hypothesis. This proves (1).
 (2) follows from (1).  
For any Dirichlet set $F$, (2) and  (2)$'$ are equivalent since $\partial \ft=-\partial F$.
(3) follows from (2) because of the positivity condition.
(3)$'$  is the Dirichlet dual of (3).
For (4) note that  if $N\ss F$ is a neighborhood of $A\in \Int F$ and $t>0$, then $tA\in \Int F$ since
$tN$ is a neighborhood of $tA$.
(5) follows from (4).
(6) follows from (3)$'$ and (5). 
(6)$'$  is the Dirichlet  dual of (6).
 \qed

\medskip\noindent
{\bf  Properties of the Class $F(X)$ for  Dirichlet-Ray Sets $F$:}
\smallskip

\item{(8)}\ \    Affine functions $u$ satisfy $\Hess_x u\in \partial F$.

\smallskip

\item{(9)}\ \    Convex functions are $F$-plurisubharmonic.

\smallskip

\item{(10)}\ \    $F$-plurisubharmonic are subaffine.

\smallskip

\item{(11)}\ \    $F$-plurisubharmonic satisfy the maximum principle.

\medskip
\noindent
{\bf Proofs:}  \  $(2) \ \Rightarrow\ (8)$, $(3) \ \Rightarrow\ (9)$, $(6) \ \Rightarrow\ (10)$.
and finally, (8) and (10) $ \ \Rightarrow\ (11)$. \qed


\vfill\eject

\noindent{\headfont  \FF. The Dirichlet Problem.}  \smallskip

In this section we state the main result, the existence and uniqueness of solutions
 of the Dirichlet problem.  We then  discuss how uniqueness follows from a local result -- 
 the Subaffine Theorem, whose proof  is postponed to Sections \GG \ and \GGG.
We conclude the section with the  proof of existence.

Given a Dirichlet set $F$ note that 
$\partial F=F\cap(\sim\Int F)=F\cap(-\ft)$, i.e., $A\in\partial F$ if and only if 
$A\in F$ and $-A\in \ft$.  Also note that $\partial \ft = -\partial F$.

\Def{\FF.1}  A function $u$ with both 
$$
u\in F(X)\qquad{\rm and}\qquad -u\in \ft(X)
$$
will be called an {\sl $F$-Dirichlet function on $X$} or an {\sl $F$-Dirichlet solution on $X$}.
\medskip
In particular, a $C^2$-function $u$ is $F$-Dirichlet if and only if 
$$
\Hess_x u\in \partial F  \qquad {\rm for\ all\ \ }  x\in X.
$$

\Theorem{\FF.2. (The Dirichlet Problem)}  {\sl  Let $\O$ be a bounded domain in $\rn$
with  smooth boundary $\bo$, and let $F$ be a Dirichlet set.
Suppose that $\bo$ is both $\Fa$ and $\oa \ft$ strictly convex. 
Then  for each $\vf\in C(\bo)$, there exists  a unique  $u\in  C(\overline\O)$ which is an $F$-Dirichlet function  on $\O$ and equals $\vf$ on $\bo$.
}

\medskip
\medskip
\noindent
{\bf Remark.}   In most interesting cases either $F\ss\ft$ or $\ft\ss F$, and then only one boundary hypothesis
is required.

\medskip
\medskip
\noindent
{\headfont Uniqueness and the  Subaffine Theorem.}   No boundary regularity is required for uniqueness,
so we replace $\overline \Omega$ by an arbitrary compact subset $K\ss\rn$

\Theorem{\FF.3. (Uniqueness)} 
      {\sl  Suppose that $F$ is a \dir set.  If  $u,v\in C(K)$ are both $F$-Dirichlet on $\Int K$
      and  $u=v$ on $\partial K$,  then $u=v$ on $K$.}
\medskip

This uniqueness theorem follows immediately from the next result.

\Theorem {\FF.4 (The Comparison Principle)}  {\sl Suppose $F$ is a  Dirichlet set
 and  that  $u, -v \in\USC(K)$.  If $$u\in F(\Int K)\and
 -v\in \ft(\Int K),$$ then}
  $$
 u\ \leq \ v \ \ \ {\rm on}\ \ \partial K\quad\Rightarrow\quad u\ \leq \ v  \ \ \ {\rm on}\ \  K
 $$
\pf
Because of the Maximum Principle (Proposition \BB.3) for subaffine functions, the 
Comparison Principle   is an immediate consequence of the next, purely local result. \qed

\Theorem{\FF.5. (The Subaffine Theorem)}  
{\sl  Assume that $F$ is a Dirichlet set.
 If $u\in F(X)$ and $v\in\ft(X)$, then $u+v \in\SA(X)$.}

\medskip 
The proof of this result    is given in Sections \GG  \ and \GGG.


\vfill\eject

\noindent
{\headfont Proof of  Existence (The Perron Solution).}  Let
$$
\cf(\vf)\ \equiv\ \left\{v\in\USC(\overline\O) : v\bigr|_\O \in F(\O)\ \ {\rm and \ \ } v\bigr|_{\bo} \leq \vf\right\}
$$
denote the Perron family for the boundary function $\vf\in C(\bo)$. 
While this family $\cf(\vf)$ may not necessarily satisfy the maximum principle,
there is a translate $\cf(\vf)+\l |x|^2$ of the family  which does.

\Lemma{\FF.6} 
 {\sl Suppose $F$ is a Dirichlet set.  Then there exists $\l>0$ with
$$
  F+\l I\ \ss\  \cpt.
$$
Hence, the maximum principle applies to $u+ \l \half |x|^2$ for all $u\in F(X)$.}
\pf
 Property (6) of Section \DD\ says that   $\ft$ is also a Dirichlet set. 
Applying Property (4) of Section \CC\ to $\ft$, pick
$ \l I\in \ft$.  This implies that $\l I +\cp\ss \ft$. Taking duals via Property (5) of Section 4 
yields $F\ss\cpt-\l I$.
\qed
\medskip

Lemma \FF.6 implies that the family 
$\cf(\vf)$  is  bounded above on $\overline \O$.
Let
$$
u(x)\ \equiv\ \sup_{v\in \cf(\vf)} v(x)
$$
denote the upper envelope of $\cf(\vf)$.

\Prop{\FF.7}  {The function $u$ belongs to $\cf(\vf)$, that is, }
$$
u\in\USC(X), \ \ \ u\bigr|_\O\in F(\O) \quad {\sl and}\quad  u\bigr|_{\bo}\leq\vf
$$
\pf 
By Property (6) in Section \DD,  $u$
 has upper semi-continuous regularization $u^*$ which satisfies
$$
u^*\bigr|_{\O }\in F(\O)
\eqno{(\FF.1)}
$$

\Lemma{\FF.8}  {\sl If $\O$ has a strictly $\oa \ft$-convex boundary, then}
$$
  u^*\bigr|_{\bo}  \ \leq    \ \vf    
  \eqno{(\FF.2)}
  $$

\medskip

 Now (\FF.1) and (\FF.2) imply that $u^*\in \cf(\vf)$.  Therefore,
$u^*\leq u$, which   is the same as 
$$
u^*\ =\ u  \quad{\rm on\ }\ \overline   \O
\eqno{(\FF.3)}
$$
This completes the proof of Proposition \FF.7 once Lemma \FF.8 is established. \qed
\medskip\noindent
{\bf Proof of Lemma \FF.8.}
By Theorem \EE.12 applied to $\ft$, there exists a global defining function 
$\rho$ which is strictly $\oa\ft$-convex on $\Ob$. 
  Pick $x_0\in \bo$.
   It follows from (\EE.4) and the affine property (2) in \S \DD\
that there exist $\e>0$ and $R>0$ so that
 $C(\rho-\e|x-x_0|^2) \in \ft(\Ob)$
if $C\geq R$.
Given $\d>0$,  pick $C>>0$ so that
$$
{\rm on \ \ } \bo:\qquad\quad      \vf + C(\rho-\e|x-x_0|^2)\ =\ \vf - C\e|x-x_0|^2\ \leq \vf(x_0)+\d.
\eqno{(\FF.4)}
$$
Then for each $v\in \cf(\vf)$
$$
w\ \equiv\ v+C(\rho-\e|x-x_0|^2)\ \in\ \SA(\O)\cap \USC(\overline\O).
$$
By the Maximum Principle we have 
$$
\sup_{\overline\O} w\ =\ \sup_{\bo} w.
$$
Now $\sup_{\bo} w \leq \vf(x_0)+\d$ since
$$
w\bigr|_{\bo}\ =\  v\bigr|_{\bo}  + C(\rho-\e|x-x_0|^2)\ \leq \vf  -  C \e  |x-x_0|^2\ \leq \vf(x_0)+\d.
$$
This proves that for all $v\in\cf(\vf)$
$$
w(x)\ =\ v(x) + C(\rho-\e|x-x_0|^2)\ \leq \vf(x_0)+\d \ \ \ \fa x\in \overline\O.
$$
Hence, the upper envelope $u$ satisfies
$$
u(x) + C(\rho-\e|x-x_0|^2)\ \leq \vf(x_0)+\d \ \ \ \fa x\in \overline\O.
$$
Therefore $u^*$ also satisfies
$$
u^*(x) + C(\rho-\e|x-x_0|^2)\ \leq \vf(x_0)+\d \ \ \ \fa x\in \overline\O.
$$
Evaluating at $x=x_0$ yields
$$
u^*(x_0) \ \leq \vf(x_0)+\d.
$$
\qed

\Lemma{\FF.9}  {\sl  If $\O$  has a  strictly $\Fa$-convex boundary, then}
$$
\liminf_{x\to x_0} u(x)\ \geq \ \vf(x_0) \qquad {\sl for\ all \ \ }  x_0\in   \bo.
$$
\pf
By Theorem \EE.12 there exists $\e>0$ so that for $C>>0$ the function $C(\rho-\e|x-x_0|^2)$ is of type $F$.
Given $\d>0$, pick $C>>0$ so that (cf. (\FF.4))
$$
{\rm on \ \ } \bo:\qquad\quad     \vf(x_0) + C(\rho-\e|x-x_0|^2)\ \leq \vf(x)+\d.
\eqno{(\FF.4)'}
$$
Set 
$$
v(x)=\vf(x_0)-\d+C(\rho-\e|x-x_0|^2)\quad{\rm on\ \ }\overline \O.
$$
 Then  $v\in \cf(\vf)$.  Consequently,
 $
 v\leq u
 $
 on $\overline\O$, and so
 $$
 \liminf_{x\to x_0} u(x) \ \geq \  \lim_{x\to x_0} v(x)\ =\    \vf(x_0)-\e. 
 $$
 \qed
 
 \Cor{\FF.10}  {\sl  If $\bo$ is both strictly $\Fa$-convex and strictly $\oa \ft$-convex, then the function
$u$ is continuous at each point of $\bo$ and $u\bigr|_{\bo}=\vf$.}
\medskip

 We now apply an argument of Walsh [W] to prove interior continuity.
 
 \Prop{\FF.11}   \qquad  $u\in C(\overline\O)$.
 
 \pf 
Let  $\O_\d \equiv\{x\in \O : \dist(x,\bo)>\d\}$  and let
 $\N\d \equiv\{x\in \Ob : \dist(x,\bo)<\d\}$.
 Suppose $\e>0$ is given.
 By the continuity of  $u$  at points of $\bo$   and the compactness of $\bo$, it follows easily
that there exists a $\d>0$ such  that:
$$
{\rm If \ \ }
|y|  \ \leq \ \d,
\ \ {\rm then\ \ } u_y \ <\ u + \e    \ \ \ {\rm on\ } \N{2\d}.
 \eqno{(\FF.5)}
 $$
where $u_y(x)\equiv u(x+y)$ is the $y$-translate of $u$ and where we
define $u$ to be $-\infty$ on $\rn-\Ob$.  We claim  that:
$$
{\rm If \ \ }
|y|  \   \leq \ \d,
\ \ {\rm then\ \ } u_y \ \leq\ u + \e    \ \ \ {\rm on\ } \Ob.
 \eqno{(\FF.6)}
 $$
This  implies,  by a change of variables, that $u\leq u_y+\e$ also holds, i.e., that: 
$$
{\rm If \ \ }
|y|  \   \leq \ \d,
\ \ {\rm then\ \ } |u_y-u| \leq \e  \ \ \ {\rm on\ } \Ob,
 $$
which completes the proof once (\FF.6) is established.

To establish (\FF.6) note first that  $u_y\in F(\O_\d)$ for each  $|y|<\d$ by Property (3) in Section 4. 
Since $u_y < u+\e$ on the collar $\N{2\d}$, one has 
$$
g_y\ \equiv\ \max\{u_y, u+\e\} \ \in\ F(\O)
$$
by Property (4) in Section 4.  Hence, $g_y - \e \in F(\O)$ by Property (2) in Section 4.
Now (\FF.5) implies that $g_y-\e =u$ on $\N{2\d}$.  Therefore,
$$
g_y-\e\ \in\ \cf(\vf)
$$
and hence $g_y-\e \leq u$ on $\Ob$.  This proves 
$$
u_y\ \leq\ g_y\ \leq\ u+\e \quad{\rm \ on\ \ } \O_\d.
$$
Combined with (\FF.5), this proves (\FF.6).\qed

 \medskip

This proves that $u\in C(\overline \O)$, $u\bigr|_\O\in F(\O)$ and $u\bigr|_{\bo}=\vf$.
To complete the proof of existence for the (DP) we show

\Lemma{\FF.12} 
$$
-u \bigr|_\O\in \ft(\O)
$$
\pf
If $-u \notin  \ft(\O)$, then since $\wt{\wt F}=F$,   Lemma \DD.6 implies that there exist $x_0\in \O$, $a$ affine, $\e>0$ and 
$A\in \Int F$ such that
$$
\eqalign{
A-u-a\ &\leq\ -\e|x-x_0|^2\quad{\rm near\ }x_0 {\rm \ \ \ and\ \ \ } \cr
&=\ \ \ 0 \ \ \qquad\qquad  \  {\rm at\ }x_0.
}
$$
Now $v=A-a+\e|x-x_0|^2$ is of type $F$. Furthermore,  $v<u$ on $\partial B_r(x_0)$ for a small $r>0$
but $v(x_0)=u(x_0)$.  Set $v'=v+\d$ with $\d>0$ small so that $v'<u$ remains true on $\partial B_r(x_0)$
but $u(x_0) < v'(x_0)$.  Then
$$
w\ =\ \cases{u \qquad \qquad \qquad {\rm on\ }\ \overline\O- B_r(x_0)  \cr
\max\{u,v'\}\ \ \ \ \ \ {\rm on\ \  }\overline{ B_r(x_0)}
}
$$
defines a function $w\in\cf(\vf)$, the Perron family for the boundary function $\vf$.
This is because the upper semicontinuity of $v'-u$ (by Proposition \DD.8) implies that
 $\{v'<u\} = \{v'-u<0\}$ is an open neighborhood of $\partial B_r(x_0)$.  However, 
$w(x_0)=v'(x_0)>u(x_0)$, contradicting the definition of $u$ as the upper envelope of $\cf(\vf)$.
\qed


 \vfill\eject

\noindent{\headfont \GG. Quasiconvex Functions.}  
\smallskip

In some sense the nicest class of $F$-\psh functions is the one where 
$F=\cp$. If $X$ is connected, then
$$
v\ \in \ \cp(X)\quad\iff \quad v\in {\rm Convex}(X)\ \ {\rm or\ \ }v\equiv -\infty.
\eqno{(\GG.1)}
$$
(See Proposition \BB.5 and its restatement as Proposition \DD.5.)
The important Property (6) in Section \DD, for ``families locally bounded above'',
can be strengthened in this case.
$$
\eqalign
{
{\rm If\ } &\cf \ {\rm is\  a\ family\ of\ convex\ functions \ which\  is \ locally\  bounded\  above,} \cr
&{\rm \qquad  then \  the\  upper\    envelope\ } v=\sup_{f\in\cf} f \ { \rm is\  also\  convex.}
}
\eqno{(\GG.2)}
$$
The point is that in this case there is no need to regularize from $v$ to $v^*$.

Another useful improvement is that:
$$
\eqalign
{
{\rm If\ } &\{ v_j\} \ {\rm is\  a\  sequence\ of \ functions \ in \ } \cp(X)\ {\rm  \ which \  converges} \cr
&{\rm \qquad  pointwise\ to\ a \ function\ } v, \ {\rm then\ }v\in\cp(X).
}
\eqno{(\GG.3)}
$$
These properties are easily established.
Another important property of convex functions is due to Alexandrov.
$$
{\rm If\ }u\ {\rm is\ a \ convex\ function, \ then}\ u\ {\rm is\ twice \ differentiable\ a.e.}
\eqno{(\GG.4)}
$$
All these properties  carry over directly to the quasi-convex case.  Their analogues, which are listed below, will be used to prove the Subaffine Theorem \FF.5.

\Def{\GG.1}  A function $u$ on $X$ is {\sl $\l$-quasiconvex} if $v=u+\l\half|x|^2$ is convex.
\medskip

Set $\cp_\l \equiv \cp-\l I$.  Then for $X$ connected, we have:
$$
u\ \in \ \cp_\l(X)\quad\iff \quad u\ \ {\rm is\ } \l{\rm-quasiconvex\ or\ \ }u\equiv -\infty.
\eqno{(\GG.1)'}
$$
$$
\eqalign
{
{\rm If\ } &\cf \ {\rm is\  a\ family\ of\ } \l{\rm-quasiconvex\ functions \ which\  is \ locally\  bounded\  above,} \cr
&{\rm \qquad  then \  the\  upper\    envelope\ } u=\sup_{f\in\cf} f \ { \rm is\  also\  } \l{\rm-quasiconvex.}
}
\eqno{(\GG.2)'}
$$
Since $\sup_{f\in\cf}(f+\l\half |x|^2) = (\sup_{f\in\cf} f)+\l\half |x|^2$, the extension  of property (\GG.3) is obvious.
$$
\eqalign
{
{\rm If\ } &\{ u_j\} \ {\rm is\  a\  sequence\ of \ functions \ in \ } \cp_{\l}(X)\ {\rm  \ which \  converges} \cr
&{\rm \qquad  pointwise\ to\ a \ function\ } u, \ {\rm then\ }u\in\cp_{\l}(X).
}
\eqno{(\GG.3)'}
$$
The extension of Alexandrov's Theorem is also obvious.
$$
{\rm If\ }u\ {\rm is\ a \ locally\ quasiconvex\ function, \ then}\ u\ {\rm is\ twice \ differentiable\ a.e.}
\eqno{(\GG.4)'}
$$
Note that if $\vf$ is smooth, then in any relatively compact subdomain there exists $\l>0$ such that
$\vf$ is $\l$-quasiconvex on the subdomain.

We will also need a final property of quasiconvex functions -- {\sl Differentiability at Maximum Points}.  Since this property is essentially vacuous for purely  convex (non-constant) functions, we include a proof.
\medskip
\noindent
{\bf (DMP):} 
\ \ 
Suppose $u$ is quasiconvex and that $x$ is a local maximum point of $u$.

\centerline
{
Then $u$ is differentiable at $x$ and $(\nabla u)(x) =0$.
}
\pf
We may assume that the maximum point is the origin and the maximum value
is zero.  Then $v(x) \equiv u(x) +\l \half |x|^2 \leq  \l \half  |x|^2$ near $x=0$ and
$v(0)=0$.  Therefore, by convexity of $v$, 
$
0=2v(0)\leq  v(x) +v(-x) \leq v(x)+\l\half|x|^2.
$
Thus, 
$$
-\l\half|x|^2\ \leq\ v(x)\ \leq\  \l\half|x|^2,
$$
which proves that $v$ is differentiable at the origin and that $(\n v)(0) = 0$.
Therefore the same conclusion holds for $u$.\qed
\medskip

Suppose now that  $v$  is a convex function.  If $v$ is differentiable
 at a point $x$, then define
$$
K(v,x)\ \equiv\ \overline{\lim_{\e\to0}}\  2 {\e^{-2}}  \sup_{|y|=1} \left\{v(x+\e y) -v(x) -\e\n v(x)\cdot y \right\}.
\eqno{(\GG.5)}
$$
Otherwise define $K(v,x)=\infty$.  If $v$ is twice differentiable at $x$, then 
$K(v,x)$ is the largest eigenvalue of $\Hess_x v$.

The next result is a key to our development.

\Theorem{\GG.2. (Slodkowski [S])}  {\sl  Suppose  $v$ is a convex function on $X$.

\noindent
If $K(v,x)\geq \L$  a.e., then $K(v,x) \geq \L$ everywhere.}
\medskip

This theorem provides a nice test for when a quasi-convex function is subaffine. Recall by
Proposition \CC.7 that $\cpt(X) = \SA(X)$ is the space of subaffine functions.

\Theorem{\GG.3}  {\sl  Suppose  $u$ is locally quasiconvex on $X$. Then}
$$
\Hess\, u\in \cpt \ \ {\rm a.e.}\quad  \Rightarrow   \quad u\in\cpt(X).
$$

\pf
Set $v(x)\equiv u(x) +\Lambda \half |x-x_0|^2$.
At a point $x$ where $u$ is twice differentiable,
$$
\Hess_x u\in\cpt \quad \iff\quad \Hess_x v \in (\cpt+\Lambda\cdot I )
 \quad \iff\quad K(v,x)\geq \Lambda.
\eqno{(\GG.6)}
$$
Thus the hypothesis $\Hess_x u\in\cpt$ a.e. is equivalent to 
$$
K(v,x)\ \geq\ \Lambda\ {\ \rm a.e. \ on\ }  X.
\eqno{(\GG.7)}
$$
By Slodkowski's Largest Eigenvalue Theorem \GG.2 this is equivalent to 
$$
K(v,x)\ \geq\ \Lambda\ {\ \rm everywhere\ on\ } X.
\eqno{(\GG.8)}
$$

We now suppose that $u\notin \cpt(X)$ and derive a contradiction.  By Lemma \BB.2
there exists $x_0\in X$, $a$ affine and $\e>0$ such that
$$
\eqalign
{
u(x)-a(x)\ &\leq \ -\e\half  |x-x_0|^2 \qquad{\rm near\ } x_0, \ \ {\rm and}  \cr
&=\ 0 \qquad\qquad\qquad\ \ \ \ {\rm at\ \ } x_0
}
\eqno{(\GG.9)}
$$
Pick $\Lambda$ so that $v(x) \equiv  u(x) +\Lambda \half |x-x_0|^2$ is convex near $x_0$.
By (\GG.9) the (DMP) implies that $u$, and hence $v$, is differentiable at $x_0$.
Thus $K(v,x_0)$ is defined by (\GG.5).
 Since
$$
\eqalign
{
v(x)-a(x)\ &\leq \  (\Lambda -\e) \half  |x-x_0|^2 \qquad{\rm near\ } x_0, \ \ {\rm and}  \cr
&=\ 0 \qquad\qquad\qquad\ \ \ \ \qquad {\rm at\ \ } x_0
}
$$
it follows  that $K(v,x_0) \leq \Lambda-\e$, a contradiction.\qed

\medskip\noindent
{\bf REMARK.}   Suppose  $v=u+\L\half|x|^2$ is convex.  Theorem \GG.3 states that: \smallskip

\qquad\qquad {\sl If  
$K(v,x)\geq \L$ a.e., then $v-\L\half|x|^2$ is subaffine

\ \ \ \qquad\qquad(or equivalently that $v$ is of type $\cpt+\L\cdot I$).}

\Cor{\GG.4. (The Subaffine Theorem for Quasiconvex Functions)}  

\noindent
{\sl   Suppose $F$ is a Dirichlet set.  If $u$ and $v$ are $\l$-quasiconvex with $u\in F(X)$ and $v\in\ft(X)$,
then $u+v\in\SA(X)$.}

\pf
By Alexandrov's Theorem, $u, v$ and $u+v$ are twice differentiable a.e., and the
a.e. hessians satisfy
$$
\Hess_x(u+v)\ =\ \Hess_x u +  \Hess_x v.
$$
By Property (7) in Section 4,  $\Hess_x u\in F$ a.e. and  $\Hess_x v\in \ft$ a.e..  Therefore,
$$
\Hess_x(u+v)\ \in\  F+\ft \ =\ \cpt \quad{\rm a.e.}.
\eqno{(\GG.10)}
$$
Since $u+v$ is $2\l$-quasiconvex, Theorem \GG.3 implies that $u+v\in \cpt(X)$, i.e.,
$u+v$ is subaffine.\qed

\medskip

Theorem \GG.3 extends from $\cpt$ to an arbitrary Dirichlet set $F$.

\Cor {\GG.5} {\sl  Suppose $u$ is   locally quasi-convex on $X$.  Then }
$$
 \Hess\, u \in F \ \ a.e.   \quad\Rightarrow\quad u\in F(X)
$$
\pf
Suppose $\Hess_x u\in F$ a.e.  Given $B\in\ft$,   $\Hess_x(u+B) 
= \Hess_xu +B\in F+\ft \ss\cpt$ a.e.  By Theorem \GG.3 this implies that $u+B\in \SA(X)$, and 
hence by Definition \DD.4 that $u\in F(X)$.\qed

Note that the converse is true as well, that is, $u\in F(X)\ \Rightarrow\  \Hess_x u \in F$ a.e..
In fact,  if $u\in F(X)$ and the second derivatives of $u$ exist at $x$, then,  as in the proof of Property  (7)  in Section 4 (the case where $u$ is $C^2$ at $x$), it follows that $\Hess_x u\in F$.


\vfill\eject

\noindent{\headfont \GGG. Sup-Convolution Approximation.}  
\smallskip
Suppose that $X$ is an open subset of $\rn$.

\Def{\GGG.1. (Sup-Convolution)}  Suppose that $u$ is a bounded function on $X$.  For each $\e>0$, define
$$
u^\e(x)\ =\ \sup_{y\in X}\left\{ u(y) - {1\over\e}|x-y|^2\right\}\qquad \forall x\in X.
\eqno{(\GGG.1)}
$$

Note that $u\leq u^\e$ on $X$.  Set $\d\equiv\sqrt{\e 2N}$ where $|u|\leq N$ on $X$, and define
$X_\d = $

\noindent  $ \{x\in X : \dist(x,\partial X) > \d\}$.
The following equivalent formulas for $u^\e$ are useful.
$$
u^\e(x)\ =\ \sup_{|x-y| \leq \d}\left\{ u(y) - {1\over\e}|x-y|^2\right\}\qquad \forall x\in X_\d.
\eqno{(\GGG.2)}
$$
\pf
If $x,y\in X$ and $|x-y|>\d$, then 
$u(y)-u(x)-{1\over\e}|x-y|^2  \leq 2N-{\d^2\over\e} =0$.
Therefore, $u(y)-{1\over\e}|x-y|^2  \leq u(x)$ if  $|x-y| > \d$.
Since $u(x) \leq u^\e(x)$, this proves that 
$$
{\sup_{|x-y|>\d, \ y\in X}}
\left\{u(y)- {1\over\e}|x-y|^2\right\} \ \leq\ u^\e(x)\quad{\rm if\ }x\in X
$$
which gives (\GGG.2).\qed

Making the change of variables $z=x-y$ in (\GGG.2) yields:
$$
u^\e(x)\ =\ \sup_{|z|\leq \d}\left\{ u(x-z) - {1\over\e}|z|^2\right\}\qquad \forall x\in X_\d.
\eqno{(\GGG.3)}
$$

\Theorem{\GGG.2. (Approximation)}  {\sl Suppose  $u\in F(X)$ with $|u(x)|\leq N$ on $X$.
Given $\e>0$, define $\d=\sqrt{2\e N}$. Then
\medskip
\item{1)} \ \ $u^\e$ decreases to $u$ as $\e\to 0$.

\medskip
\item{2)} \ \ $u^\e$  is ${1\over \e}$-quasiconvex.

\medskip
\item{3)} \ \ $u^\e \in F(X_\d)$.
}

\pf
For 1) note that $\e_1<\e_2 \iff -{1\over \e_1} <  -{1\over \e_2}$.  Now any of 
(\GGG.1), (\GGG.2) or (\GGG.3) imply that $u^\e$ is monotone decreasing as $\e\to0$.
By (\GGG.2)
$$
u^\e(x)\ \leq \sup_{|x-y|\leq\d} u(y) \quad \forall \,x\in X_\d.
$$
As noted above, $u\leq u^\e$.  Hence, 
$$
u(x)\ \leq\ u^\e(x)\ \leq\ \sup_{|x-y|\leq\d} u(y).
$$
Since $u\in \USC(X)$, the functions $\sup_{|x-y|\leq\d} u(y)$ decrease to $u(x)$. This proves 1).

To prove  2) we first note that for  $y\in X$ fixed, the function 
$u(y)-{1\over \e}|x-y|^2 +{1\over \e}|x|^2$ is affine and hence convex.  That is,
$u(y)-{1\over \e}|x-y|^2$ is a ${1\over \e}$-quasiconvex function of $x$.
Now applying (\GG.2)$'$  to (\GGG.1) (Note that (\GGG.2) does not work here) proves that 
$u^\e$ is ${1\over \e}$-quasiconvex.  

To prove  3) we make use of (\GGG.3).  Each function $u_z(x)=u(x-z)\in F(X_\d)$ if $|z|\leq \d$
by the translation property (3) in Section \DD.
Therefore, by the ``families locally bounded above'' property (6) in Section \DD,
the upper envelope $u^\e$ of the family
$$
\cf\ =\   \left\{u(x-z)-\smfrac 1 \e |z|^2 : |z|\leq\d   \right\}
$$
has upper-semicontinuous regularization in $F(X_\d)$.  However, $u^\e$ is continuous since it is quasiconvex.  Hence $u^\e$ equals its u.s.c. regularization.\qed
\medskip

The Subaffine Theorem \FF.5  follows easily from the quasi-convex case (Corollary \GG.4)
because of the Approximation Theorem  \GGG.2.

\noindent
{\bf Proof of The Subaffine Theorem \FF.5.}
The result is local, so  by upper semicontinuity we may assume $u$ and $v$ are bounded above.
We may also assume they are bounded below by replacing them with 
$u_m= \max\{u, -m\}$ and $u_m= \max\{u, -m\}$, and then taking the decreasing limit
of $u_m+v_m$ as $m\to\infty$.

We now apply Theorem \GGG.2  to $u$ and $v$ to obtain sequences 
$\{u_j\}$  and $\{v_j\}$ which are quasi-convex for each $j$ and converge
monotonically downward to $u$ and $v$ respectively as $j\to\infty$.  
By Corollary \GG.4, the sum $u_j+v_j \in \SA(X)$
for all $j$.   Since $u_j+v_j$ decreases to $u+v$, Property (5) in Section \DD, applied to the subaffine case, implies that $u+v\in \SA(X)$.
\qed


\vfill\eject

\noindent{\headfont {{\HH.}}  Topological Restrictions on Domains with Strictly $\Fa$-Convex Boundaries.}
\smallskip

In this section we show that the strict $\Fa$-convexity of $\bo$, which was assumed in our Main Theorem
\FF.2,  often places strong restrictions on the topology of $\O$.  A typical example is that
of a domain in $\bbc^n$ with pseudoconvex boundary (a Stein domain) which has the 
homotopy-type of a complex of dimension $\leq n$.  Theorem  \HH.5 greatly generalizes this fact.
For its statement we need to introduce the following ideas. 

Suppose $\rn = N\oplus W$ is an orthogonal decomposition of $\rn$.  Let $\pi_W :\Symn\to\Sym(W)$
denote restriction of quadratic functions.

\Def{\HH.1}  Suppose that $F$ is a Dirichlet-Ray set, and that $\rn = N\oplus W$.
\smallskip
\item{(1)}\ \ $W$ is  {\sl $F$-free} if $\pi_W(F)=\Sym(W)$.

\smallskip
\item{(2)}\ \ $W$ is  {\sl $F$-Morse} if there exists $A\in F$ with $\pi_W(A)<0$.

\smallskip
\item{(3)}\ \ $N$ is  {\sl $F$-strict} if $P_N\in \Int F$.

\Prop{\HH.2}  {\sl  Suppose $F$ is a Dirichlet-Ray set, and that $\rn = N\oplus W$. 
Then the following conditions are equivalent:}
\smallskip
\item{(1)}\ \ $W$ is  {\sl $F$-free},
\smallskip
\item{(2)}\ \ $W$ is  {\sl $F$-Morse},

\smallskip
\item{(3)}\ \ $N$ is  {\sl $F$-strict}
\smallskip
\noindent
The proof is given at the end of this section.

\Def{\HH.3}   The {\sl free dimension} of a \dir ray set $F$, denoted by free-dim$(F)$, is the maximal dimension of an $F$-free subspace of $\rn$.
By Proposition \HH.2, free-dim$(F)$ is also the 
maximal dimension of an $F$-Morse subspace of $\rn$.
 
 \Ex{\HH.4} Suppose $F=\cp(G)$ is   defined by a closed subset
 $G\ss G(p, \rn)$ of the Grassmannian of $p$ planes, as in (\II.10) below.  Then a subspace
 $W\ss\rn$ is $\cp(G)$-free if and only if it contains no $G$-planes, i.e.,  
 $$
\not\exists\ \x\in G \ \ \ {\rm with\ \ }\x  \ss W
 $$
(For the proof see [HL$_3$].) 
This enables one to easily calculate the free dimension in all the standard calibrated geometries. For example, when $G\ss G(2,\bbr^{2n})$ is the Grassmannian 
of complex lines in $\bbc^n=\bbr^{2n}$, the free dimension is $n$. 
This is the Stein case.
In associative geometry, the free dimension is 3, and in coassociative geometry it is 4.  When $G$ is the space of Lagrangian $n$-planes in $\bbc^n$, the free dimension is $2n-2$.

 \medskip
 
 The following theorem is the main result of this section.
It represents a surprising extension of  the Andreotti-Frankel Theorem in complex analysis 
to this  general context.

\Theorem{\HH.5}  {\sl Let $F\ss\Symn$ be a Dirichlet-Ray set with free-dim$(F) = D$. 
Suppose $\O\ss\ss\rn$ is a domain with a smooth, strictly $F$-convex boundary. Then
$\O$ has the homotopy-type of a CW-complex of dimension $\leq D$.
}

\Cor{\HH.6} {\sl Let $\O\ss\ss X$ be a   domain with a smooth, strictly $F$-convex  boundary, and 
let $D$ be the free dimension of $F$.  Then 
$$
H_k(\bo  ; \,\bbz)\ \cong \ H_k(\O  ; \,\bbz)\qquad  {\rm for\ all\ } k< n-D-1
$$
and the map $H_{n-D-1}(\bo;\, \bbz)\to H_{n-D-1}(\O;\, \bbz)$ is surjective.}

\medskip
\noindent
{\bf Proof of Theorem \HH.5}
By Theorem \EE.12 there exists a global defining function $\rho\in C^\infty(\ob)$ for $\bo$ which is
strictly  $F$-\psh on $\ob$. Set  $u = -\log(-\rho)$ and note that $u$ is a proper 
proper exhaustion function for $\O$. Direct computation shows that
$$
\Hess\, u \ = \ -{1 \over \rho} \Hess \, \rho + {1\over \rho^2} (\nabla \rho \circ \nabla \rho).
$$
Since $\Hess_x \rho\in \Int F$ and $(\nabla \rho \circ \nabla \rho)_x\in \cp$,
Property (3) shows that 
$$
\Hess_x u \in \Int F
\eqno{(\HH.1)}
$$
 at each $x\in \O$, i.e., $u$ is strictly
$F$-\psh on $\O$.
 By standard approximation 
theorems (cf. [MS]) we may assume that all critical points of $u$ are non-degenerate.
The theorem will follow from   Morse Theory if we can show that each critical point 
$x_0$ of $u$ in $X$  has index 
$\leq D$.

Suppose $x_0$ were a critical point of index $>D$.  Then there would exists 
a linear subspace $W\ss T_{x_0}\rn  = \rn$ of dimension $>D$ such that 
$$
\Hess_{x_0} u \bigr|_W\ <\ 0.
$$
However, by Proposition \HH.2  and  (\HH.1) we see that $D$ is the largest dimension of a subspace  $W$ for which this  can hold. 
 Hence, the index of $\Hess_{x_0} u \leq D$ as desired.\qed

\medskip
\noindent
{\bf Proof of Corollary \HH.6}
 This follows from the exact sequence
$$
H_{k+1}(\O, \bo  ; \,\bbz)\ \to \ H_k(\bo  ; \,\bbz)\ \to \  H_k(\O  ; \,\bbz)\ \to \ H_k(\O,\bo  ; \,\bbz),
$$
Lefschetz Duality: $H_k(\O,\bo;\,\bbz) \cong H^{n-k}(\O;\,\bbz)$, and Theorem \HH.5.\qed

\medskip
\noindent
{\bf Proof of Proposition \HH.2}
Let $\rho:\Symn\to\Sym(N)^\perp$ denote orthogonal projection onto the subspace  
$\Sym(N)^\perp$ of $\Symn$.  Also consider the conditions:
\medskip
\item{(1)$^*$}\ \ $\rho(F)\ =\ \Sym(N)^\perp$, and

\medskip
\item{(3)$^*$}\ \ $\Sym(N) \cap \Int F \ \neq\ \emptyset.$
\medskip
\noindent
The implications $(1)^* \Rightarrow (1)$ and $(3)^* \Rightarrow (3)$ are trivial.
We will prove that $(1) \Rightarrow (3)$, $(3)^* \Rightarrow (1)^*$ and $(1) \iff (2)$.
\medskip
\noindent
{\bf Proof that $(1) \Rightarrow (3)$:}
By (1) there exists $A\in F$ with $\pi_W(A)= -I_W$ where $I_W$ denotes the identity on $W$.
It suffices to show that there exist $t>0$ and $P>0$
such that $P_N = tA+P$, because by the ray property $tA\in F$, and $F+\Int\cp\ss\Int F$.
In terms of the $2\times 2$ blocking induced by $\rn= N\oplus W$, we have 
$A=\left(\matrix{a&b \cr b& -I_W}\right)$. Therefore
$$
\smfrac 1 t P \ =\ \smfrac 1 t P_N - A \ =\ \left(\matrix{\smfrac 1 t I_N -a & -b\cr -b& I_W}\right).
$$
For $t>0$ sufficiently small, we have $\smfrac 1 t P >0$ and hence $P>0$. \qed

\medskip
\noindent
{\bf Exercise}.  Show that $\left(\matrix{\smfrac 1 t I_N -a & -b\cr -b& I_W}\right) \ >\ 0$ for all sufficiently small $t>0$.

\medskip
\noindent
{\bf Proof that  $(3)^* \Rightarrow (1)^*$:}
Suppose $A\in\Sym(N)\cap \Int F$, i.e., $A=\left(\matrix{a&0 \cr 0& 0}\right)$.
Given  $B\in\Sym(N)^\perp$, i.e.,  $B=\left(\matrix{0&b \cr b& c}\right)$,
pick $\e>o$ small enough so that $A+\e B\in F$. By the ray property, ${1\over \e} A+B = 
\left(\matrix{{1\over \e} a&b \cr b& c}\right) \in F$.  Finally $\rho({1\over \e} A+B)=B$.\qed

\medskip
\noindent
{\bf Proof that  (1) $\iff$ (2) :}  Note that:
$W$ is not $F$-Morse \ $\iff$\  $\pi_W(F)\ss \cpt_W$ where $\cp_W=\{A\in\Sym(W):A\geq 0\}$
 \ $\iff$\  $\cp_W\ss\wt{\pi_W(F)}$.  
 It is easy to show that $\pi_W(F)$ satisfies the positivity condition.  
 Moreover, $\pi_W(F)\neq \emptyset$ since $F\neq \emptyset$.  Hence, $\pi_W(F)$ is either a 
 \dir set or $\pi_W(F)=\Sym(W)$.  In either case $\wt{\pi_W(F)}$ satisfies 
 the positivity condition. Hence,
 $$
 \cpt_W\ \ss\ \wt{\pi_W(F)}\quad\iff\quad 0\in \wt{\pi_W(F)}.
 $$
By definition, $0\in \wt{\pi_W(F)} \iff 0\notin \Int\, \pi_W(F)$.
Since $F$ satisfies the ray condition, so does $\pi_W(F)$.  Therefore, 
$$
\pi_W(F)=\Sym(W) \quad\iff\quad 0\in\Int \, \pi_W(F),
$$
or equivalently,
$$
0\notin\Int \, \pi_W(F)\quad\iff\quad \pi_W(F)\neq\Sym(W),
$$
i.e., $W$ is not $F$-free.\qed


\vfill\eject

\noindent{\headfont   \II.  Examples of \dir Sets. }\smallskip

Dirichlet sets  $F\ss\Symn$, to which our main existence and uniqueness 
theorem applies, are abundant, interesting and quite varied. They arise in quite different
contexts, and we have tried  to organize our presentation in that way.  There are
however, some organizational principles which illuminate the constructions.  We shall
mention these early on.

In many cases the $C^2$-solutions to the Dirichlet problem associated to $F$ satisfy 
an explicit nonlinear second-order differential equation.  When this is so, the equations
will be presented.

As mentioned in the introduction, readers are encouraged 
to look at  examples  close  to their interests and bear them in mind while reading
other parts of the paper.

\medskip
\noindent
{\bf Three Fundamental Examples.}
The most basic example of a Dirichlet set is the
set 
$$
\cp\ =\ \{A\in \Symn :   A\geq0 \},
$$
of non-negative symmetric matrices, whose \dir dual is the set
$\cpt$  of matrices with at least one non-negative eigenvalue.
These sets have analogues over $\bbc$ and $\bbh$.

  Consider the three vector spaces $\rn, \bbc^n$, and
$\bbh^n$ with scalar field $\K=\bbr, \bbc$ and $\bbh$ respectively. (In the quaternionic
case it is convenient to have the scalars $\bbh$ act  on $\bbh^n$ from the right.)
Let  $G(p, \K^n)$ denote the grassmannian of $p$-dimensional $\K$-planes in $\K^n$.
For each $\x\in G(p, \K^n)$ define the $\x$-trace of $A\in\Sym_\bbr(\K^n) =\Sym(\bbr^N)$
(with $N=n,2n$ or $4n$) by 
$$
\tr_\x A \ =\ {\rm trace}\left\{A\bigr|_\x\right\} \ =\ \langle A,P_\x\rangle
\eqno{(\II.1)}
$$
where $P_\x\in\Sym(\bbr^N)$ is orthogonal projection onto $\x$ and $\langle\cdot\,,\cdot\rangle$
is the standard inner product on $\Sym(\bbr^N)$.  Define
$$
\cp(\bbr^n) \ =\ \{A\in\Sym(\bbr^n) : \tr_\x A \geq0 \ \forall \x\in G(1,\bbr^n)\}
\eqno{(\II.2)}
$$
$$
\cp_\bbc(\bbc^n) \ =\ \{A\in\Sym_\bbr(\bbc^n) : \tr_\x A \geq0 \ \forall \x\in G(1,\bbc^n)\}
\eqno{(\II.3)}
$$
$$
\cp_\bbh(\bbh^n) \ =\ \{A\in\Sym_\bbr(\bbh^n) : \tr_\x A \geq0 \ \forall \x\in G(1,\bbh^n)\}
\eqno{(\II.4)}
$$
These are the three fundamental example of \dir sets.  
Note that they are convex cones in $\Sym(\bbr^N)$ with vertex at the origin.
Their \dir duals are given respectively by
$$
\eqalign
{
\cpt(\bbr^n) \ &=\ \{A\in\Sym(\bbr^n) : \exists\, \x\in G(1,\bbr^n) \ {\rm s.t.}\ \tr_\x A \geq0 \}     \cr
\cpt_\bbc(\bbc^n) \ &=\ \{A\in\Sym_\bbr(\bbc^n) :  \exists\, \x\in G(1,\bbc^n) \ {\rm s.t.}\  \tr_\x A \geq0 \}    \cr
\cpt_\bbh(\bbh^n) \ &=\ \{A\in\Sym_\bbr(\bbh^n) :  \exists\, \x\in G(1,\bbh^n) \ {\rm s.t.}\ \tr_\x A \geq0 \}    \cr
}
$$

Given $A\in \Sym_\bbr(\K^n)=\Sym(\bbr^N)$, consider the hermitian symmetric part of $A$.
In the complex case $\bbc^n = (\bbr^{2n},J)$ this is just 
$$
A_\bbc \ =\ \half (A-JAJ)
\eqno{(\II.5)}
$$
while in the quaternionic case $\bbh^n = (\bbr^{4n},I,J,K)$ it is
$$
A_\bbh \ =\ \smfrac 1 4(A-IAI-JAJ-KAK)
\eqno{(\II.6)}
$$
The hermitian symmetric part is $\K$-linear with $n$ eigenvalues $\l_1,...,\l_n$.
The \dir sets $\cp(\rn), \cp_\bbc(\bbc^n)$ and $\cp_\bbh(\bbh^n)$ can all be (equivalently) defined as
$$
\cp_\K(\K^n)\ =\ \{ A\in \Sym(\bbr^N) : A_\K\geq0\} \ =\ \{ A\in \Sym(\bbr^N) : \l_1\geq0,...,\l_n\geq0\} 
\eqno{(\II.7)}
$$

\vskip.3in

\noindent
{\bf The Monge-Amp\`ere Equation.}  
In all three cases there is a determinant function on $\SymN$:
$$
\det_\K A\ =\ \l_1\cdots\l_n.
\eqno{(\II.8)}
$$
Of course , if $\K=\bbr$, this is the real determinant of $A\in\Symn$,
and if $\K=\bbc$, then this  is the complex determinant of the hermitian
symmetric part of $A\in\Sym_\bbr(\bbc^n)$.  If $\K=\bbh$, then one can show that 
$\det_\bbh A$ is also a polynomial of degree $n$ (cf. [DK]).
Note that in each of these cases the boundary of the Dirichlet set 
($\partial \cp$, $\partial \cp_\bbc$ or $\partial \cp_\bbh$) is contained in the zero locus of the 
determinant  function ($\det_\bbr$,   $\det_\bbc$,   $\det_\bbh$).
Therefore, if $u$ is a $\cp_\K$-\dir function which is $C^2$, then at each point
$$
\det_\K(\Hess \, u)\ =\ 0
\eqno{(\II.9)}
$$

\vskip .3in
\noindent
{\bf The Next Tier: Other Branches of  Det(Hess u) = 0.}  Fix a positive integer $0\leq q\leq n-1$
and consider the sets
$$
\eqalign
{
P_q(\K^n) \ &=\ \{A\in\SymN  : \exists\ W\in G(n-q,\K^n)\ {\rm with\ } A\bigr|_W \in \cp_\K(W)\}  \cr 
                        &=\ \{A\in\SymN  :  A_\K\ {\rm has\ at\ least\ } n-q\ {\rm eigenvalues}\  \geq 0\}  \cr 
                        &=\ \{A\in\SymN  : \forall \ V\in G(q+1,n)\  A\bigr|_V \in \cpt_\K(V)\}.  \cr 
 }
 $$
It is easy to see that
$$
\wt P_q(\K^n) \ =\ \{A\in\SymN : \forall\ W\in G(n-q,\K^n) \ A\bigr|_W \in \cpt_\K(W)\}. 
$$
In all three cases,
$$
\wt P_q \ =\ P_{n-q-1} \qquad\qquad{\rm and \   }\qquad\qquad \cases{P_0 \ =\ \cp \cr P_{n-1}\ =\ \cpt\cr}
$$
and, therefore,  $u$ is a $P_q$-\dir function if $u\in P_q(X)$ and $-u\in P_{n-q-1}(X)$. 
Thus, if $u$ is $C^2$ and $\l_1(x)\leq \l_2(x)\leq \cdots\leq \l_n(x)$ are the eigenvalues of
$\Hess_x u$, then $u$ is $P_q$-\dir iff 
$$\l_{q+1} \equiv 0.$$
No matter what $q$ ($0\leq q\leq n-1$), one has
$$
\partial P_q\ \ss\ \{A : \det_\K A_K=0\}
$$
and, in fact, $\partial P_q$ consists of the  branch of 
$\{A : \det_K A_K =0\}$ where $\l_{q+1}=0$. In particular,  a $P_q$-\dir function which is $C^2$ satisfies the 
Monge-Amp\`ere equation (\II.9).

\vskip .3in
\noindent
{\bf  \dir Sets which are Geometrically Defined.} 
The three fundamental examples $\cp,\cp_\bbc$ and $\cp_\bbh$ are geometrically defined
by the three Grassmannians $G(1,\bbr^n), G(1,\bbc^n)$ and $G(1,\bbh^n)$ respectively.
In fact, there exists a vast array of geometrically interesting \dir sets defined in a similar fashion.
  Let $G\ss G(p,\bbr^n)$ be a closed subset of the Grassmannian of $p$-planes,
and define
$$
\cp(G) \ =\ \{A\in\Symn : \forall\, W\in G, \ \tr_W (A)\ \geq\ 0\}
\eqno{(\II.10)}
$$
This is evidently a \dir set.  It is also a convex cone with vertex at the origin. Its \dir dual is
$$
\wt \cp(G) \ =\ \{A\in\Symn : \exists \, W\in G, \ \tr_W (A)\ \geq\ 0\}
$$
 In these cases the $\cp(G)$-\psh functions have the nice
property that they are subharmonic on minimal $G$-submanifolds (those whose tangent
planes lie in $G$).  
There are many other important cases coming from calibrated geometry
and symplectic geometry.
This and other related matters are discussed in detail in [HL$_{2,4,5}$], and we briefly 
describe them next.

\vskip .3in
\noindent
{\bf  The \dir Problem in Calibrated Geometry.}
Let $\f\in \L^p\rn$ be a (constant coefficient) calibration on $\rn$, and let
$G(\f) = \{\x\in G(p,\rn): \f(\x)=1\}$ be the Grassmannian of $\f$-planes (cf. [HL$_1$]).
Then we have a geometrically defined \dir set given by (\II.10).  The attendant notions
of $\f$-\psh functions and $\f$-convexity are discussed in detail in [HL$_2$].
Our Main Theorem \FF.2   shows that on strictly $\f$-convex domains $\O\ss\rn$
the \dir problem is uniquely solvable in the class of continuous $\f$-\dir functions for 
all continuous boundary data.

We recall that this includes many interesting cases, for example, Special Lagrangian
Geometry, Associative, Coassociative and Cayley Geometries, and many others.
When a solution  to the \dir problem is $C^2$,  it is  partially $\f$-pluriharmonic, that is 
$\tr_\x \Hess\, u\geq 0$ for all $\f$-planes $\x$ and $=0$  for some $\f$-plane $\x$ 
at each point. The associated differential
equations of Monge-Amp\`ere type in these cases have not all been found.

\vfill\eject
\noindent
{\bf  The \dir Problem in Lagrangian Geometry.}
Consider $\bbc^n =(\bbr^{2n},J)$ as before, and for $A\in\Sym(\bbr^{2n})$ define
its {\sl Lagrangian component} to be
$$\eqalign{
A_{LAG} \ &=\ \smfrac t2 I + \half(A+JAJ)  \cr
&=\ \smfrac t2 I + A_{\rm skew}\cr}\qquad\qquad{\rm where\ } t=\tr_\bbr A.
$$
 The matrix $A_{\rm skew}$, called the {\sl skew-hermitian part of $A$}, anticommutes
 with $J$ and therefore has eigenvalues 
 $$
 \l_1, -\l_1,\l_2,-\l_2,...,\l_n,-\l_n
 $$
with corresponding eigenvectors of the form
$$
e_1, Je_1, e_2, Je_2,...,e_n,Je_n.
$$
Following [HL$_4$] we consider the expression
$$
{\bf M}_{LAG}(A) \ =\ \prod_{2^n\ {\rm times}} \left({t\over 2}\pm\l_1 \pm\l_2 \cdots \pm\l_n \right)
\eqno{(\II.11)}
$$
This is a polynomial in $t$ whose coefficients are symmetric functions in $\l_1^2,...,\l_n^2$.
It follows from the work of   Dadok and Katz  [DK] that ${\bf M}_{LAG}(A)$ is a polynomial 
in the coefficients of $A$.  It is, in fact, one of the factors of $\tr(D_{A_{LAG}})$ on
$\L^n\bbr^{2n}$.

We now consider the set $\LAG\ss G(n,2n)$ of Lagrangian $n$-planes in
$\bbr^{2n}=\bbc^n$.  This gives us the geometrically defined \dir set
$$
\eqalign
{
\cp(\LAG)\ &=\   \{ A\in\Symn : \forall\, \x\in \LAG,\  \tr_\x A\geq 0\}    \cr
                 &=\   \{ A\in\Symn : \smfrac t2-\l_1-\cdots-\l_n\geq 0\}    \cr
}
$$
where we assume by convention that $0\leq \l_1\leq\cdots\leq\l_n$.\ The  \dir dual is
$$
\eqalign
{
\wt\cp(\LAG)\ &=\   \{ A\in\Symn : \exists\, \x\in \LAG,\  \tr_\x A\geq 0\}    \cr
                 &=\   \{ A\in\Symn : \smfrac t2-\l_1+\cdots+\l_n\geq 0\}    \cr
}
$$

Our \dir problem on a strictly Lagrangian-convex domain $\O\ss\rn$ is uniquely solvable
for continuous boundary data and gives a  Lagrangian \psh function
$u \in C(\ob)$  which, when it is class $C^2$, satisfies the differential equation
$$
{\bf M}_{LAG}(\Hess\, u) \ =\ 0.
$$

We can now elaborate this discussion using the general principle above.
Fix a positive integer $p\leq n$ and consider the set
$$
\ISO \ =\ \{\x\in G(p, 2n) : \x \ {\rm is\  an\  isotropic\ } p\ {\rm plane}\}.
$$
(Recall that a $p$-plane $\x$ is {\sl isotropic} if $\x \perp J\x$, or equivalently, 
if $\o\bigr|_\x=0$ where $\o$ is the standard K\"ahler form.) Following the gnenral
principle we introduce the \dir sets
$$
\eqalign
{
\cp^+(\ISO)\ &=\   \{ A\in\Symn : \forall\, W\in G_\bbc(p,n),\ A\in \cp^+(\LAG)(W)\}    \cr
                   &=\   \{ A\in\Symn : \forall\, \x\in \ISO,\  \tr_\x A\geq 0\}    \cr
                 &=\   \{ A\in\Symn : \smfrac {p}{2n} t-\l_{n-p+1}-\cdots-\l_n\geq 0\}    \cr 
}
$$
and its \dir dual
$$
\eqalign
{
\wt\cp^+(\ISO)\ &=\   \{ A\in\Symn : \exists \, W\in G_\bbc(p,n),\ A\in \cp^+(\LAG)(W)\}    \cr
                   &=\   \{ A\in\Symn : \exists \, \x\in \ISO,\  \tr_\x A\geq 0\}    \cr
                 &=\   \{ A\in\Symn : \smfrac {p}{2n} t+\l_{n-p+1}+\cdots+\l_n\geq 0\}    \cr 
}
$$
Associated to this problem we have the polynomial
$$
{\bf M}_{\ISO}(A) \ =\  \prod_{|I|=p\,{\rm and\ }\pm}
 \left(\smfrac {p}{2n} t \pm\l_{i_1} \pm\cdots\pm \l_{i_p}\right)
$$
which is also a factor of $\tr(D_{A_{LAG}})$ on
$\L^n\bbr^{2n}$. As above we have that any $C^2$ function $u$ which is 
$\ISO$-Dirichlet satisfies the differential equation
$$
{\bf M}_{\ISO}(\Hess\, u) \ =\ 0
$$

\vskip .3in
\noindent
{\bf  The Geometrically $p$-Plurisubharmonic \dir Problem.} 
There is a second, more geometric, choice for the $p$-\psh functions, different from the one made at the beginning of this section. Namely, consider for $1\leq p\leq n$, the geometrically defined \dir sets
$$
\eqalign
{
\cp\left(G(p,\bbr^n)\right)\ &=\   \{ A\in\Sym(\bbr^n) :  \tr_\x A \geq0 \ \forall \x\in G(p,\bbr^n)\}    \cr
\cp\left(G(p,\bbc^n)\right)\ &=\   \{ A\in\Sym_\bbr(\bbc^n) :  \tr_\x A \geq0 \ \forall \x\in G(p,\bbc^n)\}    \cr
\cp\left(G(p,\bbh^n)\right)\ &=\   \{ A\in\Sym_\bbr(\bbh^n) :  \tr_\x A \geq0 \ \forall \x\in G(p,\bbh^n)\}    \cr
}
$$
The \dir duals are:
$$
\eqalign
{
\cpt\left(G(p,\bbr^n)\right)\ &=\   \{ A\in\Sym(\bbr^n) :   \exists\ \x\in G(p,\bbr^n) \ {\rm s.t.\ } \tr_\x A \geq0 \}    \cr
\cpt\left(G(p,\bbc^n)\right)\ &=\   \{ A\in\Sym_\bbr(\bbc^n) :  \exists\  \x\in G(p,\bbc^n) \ {\rm s.t.\ } \tr_\x A \geq0 \}    \cr
\cpt\left(G(p,\bbh^n)\right)\ &=\   \{ A\in\Sym_\bbr(\bbh^n) :   \exists\ \x\in G(p,\bbh^n) \ {\rm s.t.\ } \tr_\x A \geq0 \}    \cr
}
$$
In all three of these cases ($\bbr, \bbc$ or $\bbh$) there is a Monge-Amp\`ere polynomial
$M_p$.  First we consider the real case.
 For $A\in\Symn$,
let $D_A:\L^p\rn \to \L^p\rn$ be the extension as a derivation.  If $A$ has eigenvalues 
$\l_1,...,\l_n$ with  eigenvectors $e_1,...,e_n$, then $D_A$ has eigenvalues 
$\l_I = \l_{i_1}+\cdots+\l_{i_p}$ with  eigenvectors $e_I=e_{i_1}\wedge\cdots\wedge e_{i_p}$
where $I=(i_1,...,i_p)$ is strictly increasing.
One can prove that
$$
\eqalign
{
\cp(G(p,\rn))\ 
                 &=\   \{ A\in\Symn : D_A \geq 0\}    \cr
                 &=\   \{ A\in\Symn :  \l_I(A)\geq0\  \forall \, |I|=p\}    \cr
}
$$
and its \dir dual
$$
\eqalign
{
\cpt(G(p,\rn))\ 
                &=\   \{ A\in\Symn : D_A \ {\rm has\ at\ least\  one\ eigenvalue}\ \geq0    \}    \cr
                  &=\   \{ A\in\Symn :  \l_I(A)\geq0\   \ {\rm for\ some\ } |I|=p\}    \cr
}
$$
If $u$ is $C^2$ and $\l_1(x)\leq\l_2(x)\leq\cdots\leq\l_n(x)$ are the eigenvalues 
of $\Hess_x u$, then $u$ is $\cp(G(p,\rn))$-\dir if and only if 
$$
\l_1+\cdots+\l_p\ \equiv\ 0.
$$
Thus $C^2$-solutions to the Dirichlet problem in this case are $p$-plurisubharmonic functions
which satisfy the differential equation
$$ 
{
M_p\left(\Hess\, u \right) \ =\ \prod_{|I|=p} \l_I\ =\ 0.
}
\eqno{(\II.12)}
$$
The polynomial $M_p(A) = \prod_{|I|=p}\l_I$ is of degree $ n\choose p$ and equals $\det(D_A)$.
 For a domain $\O\ss\rn$, the Dirichlet problem for  
 $\cp(G(p,\rn))$-\dir functions can be solved uniquely  provided the boundary   is $p$-convex, i.e., 
 $$
 \tr_W\{ II_{\bo}\} \ <\ 0
 $$
for all $p$-planes tangential to $\bo$, where $II_{\bo}$ denotes the second fundamental
form of $\bo$ with respect to the outward-pointing normal.  See  [HL$_5$] for a more detailed discussion of this case, as well as a discussion of the Levi-problem in this context.

 In all three cases ($\bbr, \bbc$ or $\bbh$) 
 $$
 \eqalign
 {
 \cp(G(p,\K^n))\ &=\ \{A\in\SymN : \l_I(A_\K)\geq 0\ \forall \, |I|=p\}\qquad{\rm and}  \cr
  \cpt(G(p,\K^n))\ &=\ \{A\in\SymN : \l_I(A_\K)\geq 0\ {\rm for\ some\ } \ |I|=p\}.
 }
 $$
The polynomial $M_p$ on $\SymN$ defined by $M_p(A) = \prod_{|I|=p} \l_I(A_\K)$ of degree
$n\choose p$ provides the nonlinear differential operator exactly as in the real case.

\vskip .3in
\noindent
{\bf   The Next Tier for $\cp(G(p,\rn))$.}
Fix positive integers $p \leq q\leq n$ and consider the convex \dir set
$$
P_q(G(p,\rn))\ =\   \left\{ A\in\Symn : \exists\, W\in G(n-q,\rn),\ A\big|_W \in \cp(G(p,W))\right\}
$$
and its \dir dual
$$
\wt P_q(G(p,\rn))\ =\   \left\{ A\in\Symn : \forall\, W\in G(n-q,n),\ A\big|_W \in \wt \cp(G(p,W)) \right\}
$$
Note that
$
P_q(G(1,\rn)) \ =\ P_q$ and $\wt P_q(G(1,\rn)) \ =\ P_{n-q-1} \ =\ P_{n-q-1}(G(1,\rn)).
$
\Lemma{\II.1}  {\sl   Let $\l_1 \leq \l_2\leq \cdots\leq \l_n$ be the eigenvalues of $A_\K$.  Then}
$$\eqalign{
A\ \in\ P_q(G(p,\rn)) \quad  &\iff\quad \l_{q+1} +\cdots + \l_{q+p}\ \geq\ 0, \qquad{\rm and}  \cr
A\ \in\ \wt P_q(G(p,\rn)) \quad   &\iff\quad \l_{n-q-p+1} +\cdots + \l_{n-q}\ \geq\ 0
}
$$
The proof is straightforward. One has the following.
\Cor{\II.2}  
$$
\wt P_q(G(p,\rn)) \ =\ P_{n-q-p}(G(p,\rn)).
$$

It follows that a $C^2$-function $u$ is $P_q(G(p,\rn))$-\dir if and only if 
$$
\l_{q+1} +\cdots + \l_{q+p}\ \equiv\ 0
$$
where $\l_1 \leq \cdots\leq \l_n$ are the eigenvalues of $\Hess\, u$.
In particular, $C^2$ solutions of the \dir problem in this case are 
$p$-\psh on $q$-planes and satisfy the equation (\II.12).
In other words, they are solutions of this equation belonging to other branches
of the   locus $M_p=0$. 

This discussion holds in perfect analogy in the    complex and quaternionic cases.

\vskip .3in
\noindent
{\bf  The Next Tier Principle.}   
We have been using the following technique to generate new examples
from known ones.  Let $\cw$ be a family of subspaces of $\rn$ with a \dir set
$F_W\ss\Sym(W)$ attached to each $W\in \cw$.
Define
$$
F\ =\ \{A\in\Symn: \forall\, W\in \cw, \ A\bigr|_W\in F_W\}.
$$
One easily verifies that
$$
\ft \ =\  \{A\in\Symn: \exists\, W\in \cw, \ A\bigr|_W \in \ft_W\}
$$
\Prop{\II.3}  {\sl The sets $F$ and $\ft$ are \dir  sets.}
\medskip

This is straightforward to verify.  The examples we examine here under the heading  ``the next tier''  are of this type.  They can be elaborated to more complicated examples by repeatedly applying this principle.  For example, let $\rn=W_1\oplus\cdots\oplus W_N$ be an orthogonal decomposition and
set $\cw=\{ W_1,...,W_N\}$, $F_W=\cp(G(2,W))$.  Then 
\smallskip
\centerline{
$
F \ = \ \{A: 
\tr_\x (A)\geq 0$ for every 2-plane $\x\ss W_k$ for every $k\}$.}

\vfill\eject

\vskip .3in
\noindent
{\bf   G\aa rding Cones.}
Let $M$ be a homogeneous polynomial of degree $m$ on  $\Symn$,
and suppose  the identity $I\in\Symn$ is a {\sl hyperbolic direction for $M$} in the sense of Garding [G].  This means that for each $A\in\Symn$, the polynomial $p_A(t)=M(tI+ A)$ has exactly $m$ real  roots, and that $M(I)=1$.  Then the associated differential operator 
$$
\M(u)\ = \ M(\Hess u)
$$
will be called an  {\sl $\MA$-operator}, and the polynomial $M$ will be called an {\sl
$\MA$-polynomial.}
\medskip

G\aa rding's beautiful theory of hyperbolic polynomials states that the set
$$
\G(M)\ =\ \{A\in\Symn: M(tI+A)\neq 0\ \ {\rm for\ } t\geq0\}
\eqno{(\II.13)}
$$
is an open convex cone in $\Symn$ equal to the connected component of $\{M>0\}$ containing $I$.
The closed convex cone
$$
F_M\ =\ \{A\in\Symn: M(tI+A)\neq 0\ \ {\rm for\ } t>0\}
\eqno{(\II.14)}
$$
is the closure of $\G(M)$.  Moreover, 
$$
\partial F_M\ =\ \{A\in\Symn:  M(A) =0 {\ \rm but\ } M(tI+A)\neq 0\ \ {\rm for\ } t>0\}.
$$
 We mention that the \dir condition
 $F_M+\cp\ss F_M$ is equivalent to $\cp\ss\ F_M$ and can be stated in 
 several equivalent ways in terms of $M$:
\smallskip

1)\ \ \ $M(tI+A)\neq 0$ for all $t>0$ and $A>0$.

\smallskip

1)$'$\ \ \ $M(tI+P_e)\neq 0$ for all $t>0$ and all unit vectors $e$.

\vskip .3in
\noindent
{\bf  Symmetric Functions of Hess(u).}
A basic example of an $\MA$-polynomial on $\Symn$ is the determinant.
By the principle  above we find that each of the elementary symmetric
functions
$$
\s_{n-\ell}(A) \ =\ {1\over \ell!}{d^{\ell}\over dt^{\ell}} \, \det (A+tI)\biggr|_{t=0}
$$
is again an $\MA$-polynomial whose associated set $F_{\s_{n-\ell}}$ is
again a Dirichlet set. 

\vskip.3in

\vfill\eject
\noindent
{\bf  The Special Lagrangian Potential Equation.}
Another interesting case to which our general theory applies, comes from the polynomial
$$
Q(A)\ \equiv \ {\rm Im}\left\{  \det(I+iA)\right\}.
$$
for $A\in\Symn$.  The associated differential equation
$$
Q(\Hess\, u)\ =\ 0,
\eqno{(\II.15)}
$$ 
governs the  potential functions in the theory of Special  Lagrangian submanifolds
(cf.  [HL$_1$]).  

The locus $\{A\in\Symn : Q(A)=0\}$ has $n$  connected components, or {\sl branches},
when $n$ is even, and $n-1$ branches when $n$ is odd. Each branch is a
proper  analytic submanifold of $\Symn$.

 The Dirichlet problem for  equation (\II.15)  was treated in
[CNS] for the case where the $\Hess\, u$ is
required to lie on one of the two outermost branches.  Under this assumption, 
smooth solutions are established for smooth boundary data on appropriately 
convex domains.  In  [CNS] the authors asked whether it is possible to treat the
other branches of this equation.

We shall show that the answer is yes.  In fact we shall study the more 
general Special Lagrangian potential equation
$$
Q_{\theta}(A)\ \equiv \ {\rm Im}\left\{e^{-i\theta} \det(I+iA)\right\}.
$$
for ${\pi\over 2}< \theta\leq{\pi\over 2}$, with associated differential equation
$$
Q_{\theta}(\Hess\, u)\ =\ 0,
 \eqno{(\II.16)}
 $$

To begin we rewrite the equation $Q_\theta(A)=0$ in the form
$$
{\rm Trace}\left\{
\arctan\left({A}\right)\right\} \ =\ \theta \pm k\pi  \qquad {\rm for\ } \ k\in\bbz, |k|<{n\over 2}
 \eqno{(\II.17)}
 $$

 \Prop{\II.4} { \sl Each of the sets 
 $$
 F_c \ \equiv \   \left\{ A\in\Symn : {\rm Trace}\left\{
\arctan\left({A }\right)\right\} \geq  c\right\}
 $$
for $- {n\pi\over 2}  < c < {n\pi\over 2}$ is a \dir set with \dir dual $$\ft_c = F_{-c}.$$}

\vfill\eject

\Cor{\II.5}  {\sl Let $\O\ss\ss\rn$ be a smoothly bounded domain which is 
both $\oa F_c$ and $\oa F_{-c}$ strictly convex, with $c$ as above.  Then the \dir 
problem for continuous $F_c$-\dir functions is uniquely solvable for all continuous
boundary data on $\bo$.}

\medskip

Note that any $C^2$-function $u$, which is $F_c$-\dir, is a solution to equation
$Q_c(\Hess\, u)=0$ which lies on the branch
$$
\Hess\, u \ \in\ \partial F_c,
$$
that is, 
$$
{\rm Trace}\left\{\arctan\left({\Hess\, u }\right)\right\}  \ =\   c
$$

It is an interesting fact that the sets $F_c$ are actually starlike with respect 
to some point in their interior except for the following finite set of cases:  When $n$ is odd, 
we must  assume $\theta \neq {\pi\over 2}$,
  and for $n$ even, we assume $\theta \neq 0$. 

The special Lagrangian potential equation is, in fact, {\sl strictly elliptic} in the sense
that there is a constant $\kappa>0$ so that 
$\dist(A+P,\partial F)\geq \kappa \|P\|$  for $A\in F$ and $P\in \cp$.

We note that for $n=3$ this equation has also been treated by Yuan [Y] who established
a $C^{2,\a}$-estimate for $C^{1,1}$ viscosity solutions.

\vfill\eject


\centerline{\headfont Appendix A.}\medskip

 \centerline{\bf Dirichlet Sets Which Can Be Defined Using Fewer of the Variables in $\rn$.}
 \bigskip
 
 Suppose $F$ is a subset of $\Symn$ which {\sl can be defined using the variables 
 in a subspace $W\ss \rn$.}  That is
 $$
 F\ =\ (F\cap \Sym(W))\oplus \Sym(W)^\perp.
 $$
 Let $F_0$ denote the subset $F\cap \Sym(W)$ of $\Sym(W)$, and let $x=(x',x'')\in W\oplus W^\perp
 =\rn$ denote the variables.

 \Ex{A.1}  Let $F_0=\cp(W)$ with $W=\bbr^p$ and $p<n$. 
In this case a $C^2$-function $u$ is of type $F$ if 
$$
\sum_{j=1}^p {\partial^2 u\over \partial x_j^2}\ \geq\ 0.
$$
That is, for each fixed $x''$ the function $u(x', x'')$ of $x'$ is subharmonic.

\Remark{A.2} It is standard in the fully nonlinear theory to use the word 
``elliptic''  to include Dirichlet sets.  Then, in particular, Example A.1 is ``elliptic''.
(See [Kr] for example.)
We prefer to reserve the work elliptic for Dirichlet sets 
which can {\sl not} be defined using fewer of the variables.

 \Def{A.3}  Given a function $u(x)$ which is upper semicontinuous with values in 
 $[-\infty,\infty)$, we say that $u$ is {\sl horizontally of type $F_0$ } if for each fixed $x''$ the function
 $u_{x''}(x') = u(x',x'')$ is of type $F_0$.
 
 \medskip
 
 {\bf Elementary Properties:}
 
 \medskip
 \item{(1)} \ \ $F$ is a Dirichlet set \quad $\iff$ \quad $F_0$ is a Dirichlet set.

 \smallskip
 \item{(2)} \ \ $F$ is convex \quad $\iff$ \quad $F_0$ is  convex.

 \smallskip
 \item{(3)} \ \ $\ft = \wt{F_0}   \oplus  \Sym(W)^\perp $.

 \smallskip
 \item{(4)} \ \ $\oa F =  {\oa F_0}   \oplus  \Sym(W)^\perp$

\bigskip

If $u$ is of class $C^2$, then it is obvious that $u$ is of type $F$ if and only if $u$
is horizontally of type $F_0$.

\Theorem{A.4}  {\sl Suppose $F =  {F_0}  \oplus  \Sym(W)^\perp$ is a Dirichlet set which can be
defined using the variables in $W$.  Then $u$ is of type $F$ if and only if $u$ is horizontally of 
type $F_0$.}
 
 \Cor{A.5}  {\sl Let $F$ be as above. Then the Subaffine Theorem is true for $F$ if and only if the 
 Subaffine Theorem is true for $F_0$.}
 
 \pf Suppose $u$ is of type $F$ and $v$ is of type $\ft$.  Because of Property (3) the Theorem applies to $v$ as well as $u$.  If the Subaffine Theorem is true for $F_0$, then $u_{x''}+v_{x''}$ is a subaffine function of $x'\in \bbr^p$.  Finally, we note that if $w_{x''}(x') = w(x',x'')$ is subaffine in $x'$ (horizontally subaffine), then $w$ is subaffine in $x=(x',x'')$.  Set $B=B'\times B'' \ss \bbr^p\times \bbr^{n-p}$.
   If $w\leq a$ on $\partial B$, then $s_{x''}\leq a_{x''}$ on $\partial B'\times B''\ss\partial B$.
   Hence, $w_{x''}\leq a_{x''}$ on $B'\times B''$. \qed
 
 \medskip
 \noindent
 {\bf Proof of Theorem A.4.} 
 Suppose that $u$ is horizontally of type $F_0$.  To show that $u$ is of type $F$ we must show
 that $u+B$ is subaffine for each $B\in \ft$.  Since $\ft = \ft_0\oplus \Sym(W)$, we have 
 $B(x',x'')=b(x')+a(x',x'')$ where $b\in \ft_0$ and $a$ is affine.  By hypothesis, $u_{s''}(x')+b(x')$ is subaffine in $x'$.
 Since $a$ is an affine function, $u+B$ is horizontally subaffine.  As noted in the proof of the Corollary,
 this implies that $u+B$ is subaffine.
 
 Suppose $u(x',x'')$ is not of type $F_0$ for some fixed $x_0''$.  
 We may assume $x_0''=0$ and that there exist $\e>0$, $x_0'=0$, and $b\in \ft_0$ such that 
 $$
 \eqalign{
 u(x',0)+b(x')     &\leq\ -\e|x'|^2\quad{\rm near \ } x'=0\cr
 &=\ 0\qquad \qquad{\rm at\ } x'=0
 }
 \eqno{(1)}
 $$
 after modifying by an affine function of $x'$ and translating so that $x_0'=0$.
 
 Consider $B(x',x'')=b(x')-\L|x''|^2$ with $\L>>0$.  By (1) and upper semicontinuity,
 $$
 u(x',x'')+B(x',x'')<0 \ \ \ {\rm on}\ \ \ |x'|=r', |x''|\leq r''
 $$
 for some $r''>0$ small.  Pick $\L$ large enough so that 
 $$
 u(x',x'')+B(x',x'')<0 \ \ \ {\rm on}\ \ \ |x'|=\leq r', |x''| = r''.
 $$
 
Since $u+B$ equals zero at $x=0$, it is not subaffine and hence $u$ is not of type $F$.
\qed

 \vfill\eject

 \def\CD{B}

\centerline{\headfont Appendix \CD.}\medskip

 \centerline{\headfont A Distributional Definition of Type F for Convex \dir Sets F.}
 \bigskip

Suppose $H$ is a closed half space in $\Symn$.  Then $H$ can be defined  by
$$
H\ =\ \{B\in \Symn : \langle A,B\rangle\geq c\}
\eqno{(\CD.1)}
$$
for some non-zero $A\in\Symn$ and some $c\in\bbr$.
Note that
$$
H\ \ {\rm is\ a\ \dir \ set\ } \quad\iff\quad A\in \cp
\eqno{(\CD.2)}
$$
since, with $B_0\in\partial H$, one has
$
c\leq \langle A,B_0+P\rangle  =  \langle A,B_0\rangle +\langle A,P\rangle  = c+\langle A,P\rangle
$
for all $P\geq 0$ if and only if $0\leq \langle A,P\rangle$ for all $P\geq 0$.
 
 Similarly, one can prove that:

\Lemma {\CD.1} {\sl If $F$ is a Dirichlet set contained in a closed half-space $H$, 
then $H$ is a Dirichlet set.}
 
 \medskip
 As a consequence of this  Lemma  we can state the 
Hahn-Banach Theorem in the context of Dirichlet sets as follows.
 
 \Cor{\CD.2} {\sl  $F$ is a convex Dirichlet set if and only if $F=\bigcap_\a H_\a$ over 
 all Dirichlet supporting half-spaces $H_\a$ for $F$.}\medskip
 
 The Dirichlet dual statement is also true.

\Lemma {\CD.3}  {\sl  If $F$ is a convex Dirichlet set, then $\ft = \bigcup_\a \wt H_\a$ over 
 all Dirichlet supporting half-spaces $H_\a$ for $F$.}

\pf     If $F\ss H_\a$, then $\wt H_\a \ss \ft$, so 
we only need to show that $\ft\ss  \bigcup_\a \wt H_\a$.
Suppose  $B\in \ft$, i.e., $-B\notin\Int F$.  We claim there exists   $H_\a$ with $-B\notin \Int H_\a$,
i.e., with $B\in \wt H_\a$.  There are two cases. If $-B\notin F$ we can pick $H_\a $ with $-B\notin H_\a$. If $-B \in\partial F$, then we can pick a supporting hyperplane $H_\a$ for $F$ at $-B$ by the Hahn-Banach Theorem, so that $-B\in\partial H_\a$.\qed

\Cor{\CD.4}  {\sl A function $u$ is of type $F$ if and only if $u$ is of type $H_\a$ on $X$
for all \dir supporting half-spaces $H_\a$ for $F$.}

\pf
Since  $F\ss H_\a$, type $F$ implies type $H_\a$.  
Conversely, if $u$ is type $H_\a$ for all supporting half-spaces $H_\a$, then $u+B$ is subaffine for all $B\in { \wt H}_\a$ and hence by Lemma \CD.3, $u$ is of type $F$.
\qed

\Prop{\CD.5}  {\sl Suppose $H$ is a half-space through the origin, defined by
\smallskip
\centerline
{
$H\ =\ \{B : \langle A, B\rangle \geq 0\}$ \quad with $A$ positive definite.
}
\smallskip\noindent
  Then the following are equivalent:
\smallskip

\item {1)}  \ \ $u$ is of type $H$,
\smallskip

\item {2)}  \ \ $u$ is  sub-$\D_A$-harmonic,
\smallskip

\item {3)}  \ \ $u$ is $L^1_{\rm loc}$ and $\D_A u \geq0$ (or $u\equiv -\infty$)

\smallskip
\noindent
where $\D_A u = \sum a_{ij} u_{ij}$.}

\pf  The equivalence of 2) and 3)  is standard.  Note that $H$ is self dual,
 i.e. $\wt H=H$.  Suppose $u$ is of type $H$.  
 Given a $\D_A$-harmonic function $h$ with $u\leq h$ on $\partial B$, we have  $u\leq h$ on $B$ because $-h$ is of type $\wt H$ which implies that $u-h$ is subaffine.
 
 Conversely, suppose $u$ is  sub-$\D_A$-harmonic.  Let $v$ be a $C^2$-function of type 
 $\wt H=H$.  We must show that $u+v\leq a$
 on $\partial B$ implies $u+v\leq a$ on $B$  for any affine function $a$ and any ball $B$.
 Replace $v$ by $v-a$ and $a$ by $0$. Let $h$ denote the $\D_A$-harmonic function with the
 same boundary values as $v$ on $\partial B$.  Now $u+h\leq0$ on $\partial B$ implies $u+h\leq 0$ on $B$  since $u$ is sub-$\D_A$-harmonic, but $v=h$ on $\partial B$
 implies $v\leq h$  on $B$ since, as  we have shown above, $v$ is sub-$\D_A$-harmonic.
 \qed

\Cor{\CD.6}  {\sl Suppose $H$ is a Dirichlet half-space defined by (\CD.1) with $A>0$.
Pick $B_0\in\partial H$. Then $u$ is of type $H$ if and only if 
$u\in\lloc$ and $\D_A(u-B_0)\geq 0$, i.e., $u-B_0$ is $\D_A$-subharmonic.}

 \Lemma{\CD.7}  {\sl A convex \ds $F$ cannot be defined using fewer of the variables in $\rn$
 if and only if each $A\in \Int \cp_+(F)$ is positive definite, where $\cp_+(F)$ denotes the closure 
 of the cone of directions defining the supporting half-spaces for $F$.}
 
 \pf   See Corollary C.4 in Appendix C of [HL$_3$].
 
 \medskip
 
 Combining Corollary \CD.4 and Lemma \CD.7 we have:
 
 \Theorem {\CD.8} {\sl Suppose $F$ is a convex \ds which can not be defined using fewer 
 of the variables in $\rn$. For each supporting half-space 
 $\{B: \langle A_\a, B\rangle \geq c\}$ pick $B^\a_0\in \partial H_\a$.
 Then $u$ is of type $F$ if and only if $u-B^\a_0$ is $\D_{A_\a}$-subharmonic
 for each $A_\a$.}
 
 
  \Remark{\CD.9} This theorem can be extended to the case where $F$ can be defined using
  fewer of the variables by applying Theorem A.4.  Moreover, one can deduce from this extension
 that for a Dirichlet set $F$
 $$
 F\ \ {\rm is\ convex} \quad\Rightarrow\quad  F(X)\ \ {\rm is\ convex} 
 $$


 \vfill\eject

\def\bbc{{\bf C}}

\centerline{\bf References}

\vskip .2in

\noindent
\item{[Al]}   S. Alesker,  {\sl  Quaternionic Monge-Amp\`ere equations}, 
J. Geom. Anal., {\bf 13} (2003),  205-238.
 ArXiv:math.CV/0208805.  

\smallskip

\noindent
\item{[Alex$_1$]}  A. D. Alexandrov, {\sl Almost everywhere existence of the second differential of a convex function and properties of convex surfaces connected with it (in Russian)}, 
Lenningrad State Univ. Ann. Math.  {\bf 37}   (1939),   3-35.

 \smallskip
%
%
%

\noindent
\item{[Alex$_2$]}    A. D. Alexandrov,  {\sl  The Dirichlet problem for the equation Det$\| z_{i,j}\| = \psi(z_1,...,z_n,x_1,...,x_n)$}, I. Vestnik, Leningrad Univ. {\bf 13} No. 1, (1958), 5-24.

\smallskip

\noindent
\item{[BT]}   E. Bedford and B. A. Taylor,  {The Dirichlet problem for a complex Monge-Amp\`ere equation}, 
Inventiones Math.{\bf 37} (1976), no.1, 1-44.

\smallskip

 \item{[Br]}  H. J. Bremermann,
    {\sl  On a generalized Dirichlet problem for plurisubharmonic functions and pseudo-convex domains.
    Characterization of \v Silov boundaries},
          Trans. A. M. S.  {\bf 91}  (1959), 246-276.
\medskip


\noindent
 \item{[CNS]}   L. Caffarelli, L. Nirenberg and J. Spruck,  {\sl
The Dirichlet problem for nonlinear second order elliptic equations, III: 
Functions of the eigenvalues of the Hessian},  Acta Math.
  {\bf 155} (1985),   261-301.

 \smallskip

%
%
%

 \smallskip

\noindent
\item{[CIL]}   M. G. Crandall, H. Ishii and P. L. Lions {\sl
User's guide to viscosity solutions of second order partial differential equations},  
Bull. Amer. Math. Soc. (N. S.) {\bf 27} (1992), 1-67.

 \smallskip

\noindent
\item{[DK]}   J. Dadok and V. Katz,   {\sl Polar representations}, 
J. Algebra {\bf 92} (1985) no. 2,
504-524.

\smallskip

\noindent
\item{ [G]}   L. G\aa rding, {\sl  An inequality for hyperbolic polynomials},
 J.  Math.  Mech. {\bf 8}   no. 2 (1959),   957-965.

 \smallskip

 \noindent 
\item {[HL$_1$]}   F. R. Harvey and H. B. Lawson, Jr,  {\sl Calibrated geometries}, Acta Mathematica 
{\bf 148} (1982), 47-157.

 \smallskip

%
%


\item {[HL$_{2}$]} F. R. Harvey and H. B. Lawson, Jr., 
{\sl  An introduction to potential theory in calibrated geometry},  
ArXiv:math.DG/0710.3920.     
\smallskip

\item {[HL$_{3}$]} F. R. Harvey and H. B. Lawson, Jr., {\sl  Plurisubharmonicity in a general geometric context},  ArXiv:math.DG/0710.3921.

\smallskip

\item {[HL$_{4}$]} F. R. Harvey and H. B. Lawson, Jr., {\sl  Lagrangian plurisubharmonicity and convexity},  Stony Brook Preprint (2007).

\smallskip

\item {[HL$_{5}$]} F. R. Harvey and H. B. Lawson, Jr., {\sl  Foundations of  $p$-convexity 
and $p$-plurisubharmonicity in riemannian geometry},  Stony Brook Preprint (2006).

\smallskip

   \noindent
\item{[HM]}    L. R. Hunt and J. J. Murray,    {\sl  $q$-plurisubharmonic functions 
and a generalized Dirichlet problem},    Michigan Math. J.,
 {\bf  25}  (1978),  299-316. 

\smallskip

   \noindent
\item{[I]}    H. Ishii,    {\sl  On uniqueness and existence of viscosity solutions of fully nonlinear second-order elliptic pde's},    Comm. Pure and App. Math. {\bf 42} (1989), 14-45.

\smallskip

   \noindent
\item{[IL]}    H. Ishii and P. L. Lions,    {\sl   Viscosity solutions of fully nonlinear second-order
elliptic partial differential  equations},    J. Diff. Eq.  {\bf 83} (1990), 26-78.

\smallskip

   \noindent
\item{[J]}    R. Jensen,    {\sl  Uniqueness criteria for viscosity solutions of fully nonlinear 
elliptic partial differential  equations},    Indiana Univ. Math. J. {\bf 38}  (1989),   629-667.

\smallskip

   \noindent
\item{[Kr]}    N. V. Krylov,    {\sl  On the general notion of fully nonlinear second-order elliptic equations},    Trans. Amer. Math. Soc. (3)
 {\bf  347}  (1979), 30-34.

\smallskip

\noindent
\item{[MS]}   J. Milnor and J. Stasheff, 
{    Morse Theory},    
Princeton University Press, Princeton, 19??.

\smallskip

\item {[RT]} J. B. Rauch and B. A. Taylor, {\sl  The Dirichlet problem for the 
multidimensional Monge-Amp\`ere equation},
Rocky Mountain J. Math. {\bf 7}    (1977), 345-364.

\smallskip

\item {[S]}  Z. Slodkowski, {\sl  The Bremermann-Dirichlet problem for $q$-plurisubharmonic functions},
Ann. Scuola Norm. Sup. Pisa Cl. Sci. (4)  {\bf 11}    (1984),  303-326.

\smallskip

\item {[So$_1$]}  P. Soravia, {\sl  On nonlinear convolution and uniqueness of viscosity solutions},
Analysis  {\bf 20}    (2000),  373-386.

\smallskip

\item {[So$_1$]}  P. Soravia, {\sl  Uniqueness results for fully nonlinear degenerate elliptic equations with discontinuous coefficients},
Comm. in Pure and Applied Analysis  {\bf 5}    (2006),  213-240.

\smallskip

\item {[W]}   J.  B. Walsh,  {\sl Continuity of envelopes of plurisubharmonic functions},
 J. Math. Mech. 
{\bf 18}  (1968-69),   143-148.

\item {[Y]}  Yu Yuan,  {\sl A priori estimates for solutions of fully nonlinear special lagrangian equations},
 Ann Inst. Henri Pioncar\'e non li\'eaire
{\bf 18}  (2001),   261-270.

\end


\centerline{\bf Appendix C.}\medskip

 \centerline{\bf Dirichlet Invariance.}
 \bigskip

If   Theorem \FF.2 holds  for a given Dirichlet set $F$, we say that:
\smallskip
\centerline{\sl The (DP) is uniquely solvable for $F$}\smallskip

\Lemma{C.1}  {\sl The (DP) is uniquely solvable for $F$ if and only if 
 the (DP) is  uniquely solvable for $F+A$ where $A\in\Symn$ is a fixed quadratic function.}

\pf Define $F'=F+A$. Then 
$$
u+A \in F'(X)\quad \iff\quad u\in F(X)
\eqno{(C.1)}
 $$
 This follows from Definition \CC.6 and the fact that $\wt{F'} = \ft+A$. 
 Setting $u'=u+A$ and $\vf'= \vf+A$ gives
 $$
u' \in F'(X)\quad \iff\quad u\in F(X)
\eqno{(C.1)'}
 $$
 and 
 $$
-u' \in \wt{F'}(X)\quad \iff\quad -u\in \ft(X).
\eqno{(C.2)}
 $$
 Applying (C.1) to $\wt{F'}=\ft-A$
 yields (C.2) since $-u'=-u-A$.\qed

 \Note{C.2} The constant function $u=c$ may not be a solution when $\vf=c$, but if $u$ is a solution for $\vf$, then $u+c$ is a solution for $\vf+c$.

 \Lemma{C.3}  {\sl The (DP) is uniquely solvable for $F$ if and only if 
 The (DP) is uniquely solvable for $\ft$.}
 
 \pf
 Suppose the (DP) is   solvable for $F$ on $\O$.  That is, suppose that for all
 $\vf\in C(\partial \O)$, there exists a unique $u\in C(\ob)$ such that:
 $$\eqalign
 {
 &{\rm On\ }\O, \ \  \ u\ \ {\rm is \ of\ type\ } F \ \ {\rm and\ } -u\ \ {\rm is \ of\ type\ } \ft
 \cr
 &{\rm On\ }\bo, \ u=\vf
 \cr
 }
 $$
 Now given $\psi\in C(\bo)$, set $\vf=-\psi$ and solve the $F$ Dirichlet Problem.  Define 
 $v= - u$ on $\ob$.  Then $v\in C(\ob)$ and we have
 $$\eqalign
 {
 &{\rm On\ }\O, \ \  \ v\ \ {\rm is \ of\ type\ } \ft \ \ {\rm and\ } -v\ \ {\rm is \ of\ type\ } \wt{\ft} = F
 \cr
 &{\rm On\ }\bo, \ v= \psi
 \cr
 }
 $$
 The uniqueness of $u$ implies the uniqueness of $v$ by reversing the process.\qed

 We next investigate  the equivariance of the Dirichlet Problem under linear automorphisms
 of $\rn$.
 
 \Lemma{C.4}  {\sl  Suppose $L:\rn\to\rn$ is a linear change of coordinates sending 
 $\O'$ to $\O$.  Then the (DP) for  $F$ is uniquely solvable on $\O$ if and only if 
 the (DP) for $F'\equiv L^tFL$ is uniquely solvable on $\O'$.}
 
 \pf  Note that  $F$ is a Dirichlet set if and only if  $F'$ is also a Dirichlet set, and that 
 $B\in\ft \iff L^tBL \in L^t\ft L= \ft'$.  If $u+B\in \SA(X)$ for all $B\in \ft$, then $u\circ L + B\circ L
 \cong L^t u L+L^tBL \in \SA(X')$, with $X'=L^{-1}(X)$, for all $B\in \ft$.  This proves that
 $u\in F(X) \iff u\circ L \in F'(X')$. The assertion now follows immediately. \qed

\end